\documentclass{elsart}

\input xy
\xyoption{all} \CompileMatrices

\usepackage{amssymb}

\newcommand{\demo}{ {\it   Proof. }}

\usepackage{theorem}
\theoremstyle{change}
{\theorembodyfont{\rmfamily}%
\newtheorem{Dfn}{Definitions}[section]

\newtheorem{Emp}[Dfn]{Example}
\newtheorem{Hip}[Dfn]{Hypothesis}
\newtheorem{Not}[Dfn]{Notation}
}
{\theorembodyfont{\itshape}
\newtheorem{Thm}[Dfn]{Theorem}

\newtheorem{Prop}[Dfn]{Proposition}
\newtheorem{Lem}[Dfn]{Lemma}
\newtheorem{com}[Dfn]{}
}

\newenvironment{apartats}{%
  \begin{enumerate}
}{%
\end{enumerate}}

\newcommand{\fix}[1]{\mbox{\rm Fix } #1}
\newcommand{\fixp}[1]{\mbox{\rm Fix } (#1)}

\begin{document}

\begin{frontmatter}



\title{A description of auto-fixed subgroups in a free group}


\author[AM]{A.\ Martino}\ead{AMartino@crm.es}, \author[EV]{E.\ Ventura}\ead{enric.ventura@upc.es}

\address[AM]{Dept.\ of Math.,\ Southampton University,\ Southampton,\ U.\ K.} \address[EV]{Dept. Mat. Apl. III,\ Univ. Polit\`ecnica
Catalunya,\ Barcelona,\ Spain \\ and \\ Dept. of Math., City College of New York, CUNY, USA}

\begin{abstract}
Let $F$ be a finitely generated free group. By using Best\-vi\-na-Handel theory, as well as some further improvements, the eigengroups
of a given automorphism of $F$ (and its fixed subgroup among them) are globally analyzed and described. In particular, an explicit
description of all subgroups of $F$ which occur as the fixed subgroup of some automorphism is given.
\end{abstract}

\end{frontmatter}

\section{Introduction}

For all the paper, let $F$ be a finitely generated free group.

The {\it rank} of $F$, denoted $r(F)$, is the cardinality of a free generating set of $F$. The {\it reduced rank} of $F$, denoted
$\tilde{r}(F)$, is $\max \{ r(F)-1, 0\}$, that is, one less than the rank, except for the trivial group where the reduced rank
coincides with the rank, which is zero. It is well known that every subgroup of a free group is free, and so it has its own rank and
reduced rank. However, they are not in general bounded above by those of $F$.

As usual, $Aut(F)$ denotes the {\it automorphism group} of $F$, $Inn(F)$ is the subgroup of {\it inner automorphisms}, and
$Out(F)=Aut(F)/Inn(F)$ is the {\it outer automorphism group} of $F$. So, an outer automorphism is a coset of $Inn(F)$ and is to be
thought of as a set of automorphisms obtained by composing a given one with all possible inner automorphisms.

Let $\phi \colon F\rightarrow F$ be an endomorphism of $F$. We will denote $\phi$ as acting on the right of the argument, $x\mapsto
(x)\phi$ (and parentheses will be omitted if there is no risk of confusion). A subgroup $H\leq F$ is called {\it $\phi$-invariant}
when $H\phi =H$, setwise. In this case, the restriction of $\phi$ to $H$ will be denoted $\phi_{_H} \colon H\rightarrow H$, and it is
an endomorphism of $H$.

Except when $F$ has rank 1, $Inn(F)$ is isomorphic to $F$. For any $y \in F$, we shall write $\gamma_{y}$ to denote the inner
automorphism of right conjugation by $y$ (denoted by exponential notation). Thus $\gamma_y \colon F \rightarrow F$, $x\mapsto
x\gamma_y=y^{-1}xy=x^y$. Similarly, for any subgroup $H\leq F$, we write $H^y=y^{-1} H y$ for its right conjugate by $y$. We will
denote the conjugacy class of $H$ in $F$ by $[[H]]$. Since $r(H)=r(H^y)$ for every $y\in F$, the {\it rank} of a conjugacy class of
subgroups, $r([[H]])$, is well defined.

The {\it fixed subgroup} of an endomorphism $\phi$ of $F$, denoted $\fix{\phi}$, is the subgroup of elements in $F$ fixed by $\phi$:
 $$
\fix{\phi}=\{ x \in F : x\phi =x \}.
 $$
For example, if $y$ is not a proper power (in particular $y\neq 1$) then, for every integer $r\neq 0$, $\fix{\gamma_{y^r}} =\langle
y\rangle$, the centralizer of $y^r$ in $F$. Note that, if $H$ is a $\phi$-invariant subgroup of $F$, then
$\fix{\phi_{_H}}=H\cap\fix{\phi}$. Following~\cite{MV1}, a subgroup $H\leq F$ is called {\it 1-auto-fixed} (resp. {\it 1-endo-fixed})
when there exists an automorphism (resp. endomorphism) $\phi$ of $F$ such that $H=\fix{\phi}$.

Following~\cite{DV}, the {\it eigengroup} of $\phi$ with {\it eigenvalue} $y\in F$ is the maximal subgroup of $F$ where $\phi$ acts as
left conjugation by $y$,
 $$
\{ x\in F : x\phi =yxy^{-1}\}=\{ x\in F : y^{-1}(x\phi) y=x\}=\fix{\phi \gamma_y}.
 $$
So, the eigengroups of $\phi$ are the fixed subgroups of the automorphisms in the outer automorphism $\Phi$ containing $\phi$.

In our view, the three main properties known about 1-auto-fixed subgroups are the following ones (and note that (iii) implies (ii)):
\begin{apartats}
\item It is easy to see that every 1-endo-fixed subgroup $H$ is
{\it pure}, i.e. $x^r \in H$ implies $x\in H$. \item The main result in~\cite{BH} states that 1-auto-fixed subgroups of $F$ have rank
at most $n=r(F)$, proving the Scott conjecture. The same result was proved for 1-endo-fixed subgroups in~\cite{IT}. \item
In~\cite{DV}, the previous result was generalized to say that every 1-mono-fixed subgroup $H$ of $F$ is {\it $F$-inert}, i.e. $r(H\cap
K)\leq r(K)$ for every $K\leq F$. Inertia for 1-endo-fixed subgroups is an open problem (see~\cite{B} and~\cite{V} for related
results).
\end{apartats}

Inertia is a quite restrictive condition, since the rank of intersections of subgroups of $F$ in general can behave like the order of
the product of the ranks. However, properties (i) and (iii) are not enough to characterize 1-auto-fixed subgroups. For example, the
subgroup $H=\langle a^{-1}b^{-1}ab, a^{-1}c^{-1}ac \rangle$ of the free group $F$ on $\{a,b,c\}$ is pure and $F$-inert, but it is not
1-endo-fixed (see~\cite{MV1}).

The goal of this paper is to provide a description for the fixed subgroups of automorphisms of free groups of finite rank. Although we
do not obtain a complete characterization of all 1-auto-fixed subgroups, we hope to provide an explicit enough answer to the question
``What subgroups $S$ of $F$ can be of the form $\fix{\beta}$~?" asked by J. Stallings in his paper~\cite{St}, section~1 problem~P2, in
1987.

Let $\{ a_1, \ldots ,a_n\}$ be a basis for $F$. It is easy to see that the automorphism $\phi$ of $F$ defined by $a_i \phi =a_i^{-1}$
for $i=1,\ldots ,n$, has trivial fixed subgroup. So, the trivial subgroup is 1-auto-fixed. Also, a cyclic subgroup $H=\langle
y\rangle$ of $F$ is 1-auto-fixed if, and only if, it is pure and, in this case, $H=\fix{\gamma_y}$. So, the interesting cases begin
with subgroups of rank 2.

The maximal rank case was completely settled by Collins and Turner. In~\cite{CT2}, these authors gave a complete description of the
1-auto-fixed subgroups $H\leq F$ with $r(H)=n$. In the current paper we generalise this, finding a similar description which applies
to all 1-auto-fixed subgroups without restriction. For later use, we reformulate part of Collins-Turner result here.

\begin{Thm}[Collins--Turner,~\cite{CT2}]\label{maxrank}
Let $F$ be a non-cyclic finitely generated free group and $\phi\in Aut(F)$ such that $r(\fix{\phi})=r(F)$. Then, there is a
non-trivial free factorisation $F=H*\langle y\rangle$, where $H$ is a $\phi$-invariant subgroup and one of the following holds:
\begin{apartats}
\item $y\phi =y$ and $\fix{\phi}=(H\cap \fix{\phi})*\langle
y\rangle$, \item $y\phi =h^r y$ and $\fix{\phi}=(H\cap \fix{\phi})*\langle y^{-1}hy \rangle$, for some $1\neq h\in H\cap \fix{\phi}$
not a proper power, and some integer $r\neq 0$.
\end{apartats}
\end{Thm}

Our main result describing 1-auto-fixed subgroups in general, is inspired by the following three basic constructions of automorphisms
and their corresponding fixed subgroups, from simpler automorphisms.

\begin{com}{\bf The basic constructions}\label{construct}
\upshape Let $H$ and $K$ be finitely generated free groups.
\begin{apartats}
\item Let $\{ a_1, \ldots ,a_m\}$, $m<n$, be a basis for $H$ and
let $\varphi \in Aut(H)$. By adding $n-m\geq 1$ extra generators $\{ a_{m+1},\ldots ,a_n\}$, we obtain a bigger free group $F$, and
$\varphi$ can be extended to an automorphism $\phi \in Aut(F)$ by setting $a_i\phi=a_i\varphi$ if $1\leq i\leq m$, and $a_i\phi =w_i$
if $m+1\leq i\leq n$, with appropriate words $w_i\in F$. Then, denoting the restriction of $\phi$ to $H$ by $\phi_{_H}$, one has
$\phi_{_{H}}=\varphi$ and it is possible to choose the $w_i$ such that $\fix{\phi}=\fix{\varphi}$ (take, for example, $w_i=a_i^{-1}$,
$m+1\leq i\leq n$). Hence, $\fix{\phi}$ is contained in $H$, which is a proper $\phi$-invariant free factor of $F$.
\item Now, let $\phi_1 \in Aut(H)$ and $\phi_2 \in Aut(K)$.
Consider $F=H*K$ and the automorphism $\phi_1 *\phi_2 \in Aut(F)$, which extends $\phi_1$ and $\phi_2$. Then, both $H$ and $K$ are
$(\phi_1*\phi_2)$-invariant, and $\fixp{\phi_1 *\phi_2}=\fix{\phi_1}*\fix{\phi_2}$.
\item Finally, let $\{ a_1, \ldots ,a_{n-1}\}$ be a basis for $H$,
let $\varphi \in Aut(H)$ and suppose that $h\varphi =h'hh'^{-1}$ for some $1\neq h,h'\in H$ with $h$ not a proper power. By adding a
new free generator $y$, we obtain a bigger free group $F$, and $\varphi$ can be extended to an automorphism $\phi \in Aut(F)$ by
setting $a_i\phi=a_i\varphi$ for $i=1,\ldots,n-1$, and $y\phi =h'h^r y$. Then, $H$ is $\phi$-invariant, $\phi_{_{H}}=\varphi$ and, for
all but finitely many choices of the integer $r$, the fixed subgroup of $\phi$ is precisely $\fix{\phi}=\fix{\varphi}*\langle y^{-1}hy
\rangle$ (see the appendix for a proof of this fact).
\end{apartats}
\end{com}

The following theorem says that the three methods above are enough to construct all non-cyclic 1-auto-fixed subgroups from simpler
ones. In other words, any automorphism of $F$ with non-cyclic fixed subgroup can be realised using one of the previous three methods
and, hence, its fixed subgroup has the corresponding decomposition.

\begin{Thm}\label{mainconnex}
Let $\phi$ be an automorphism of a finitely generated free group $F$. Then, either $\fix{\phi}$ is cyclic or there exists a
non-trivial free factorisation $F=H*K$ such that $H$ is $\phi$-invariant and one of the following holds:
\begin{apartats}
\item $\fix{\phi}\leq H$, \item $K$ is also $\phi$-invariant and
$\fix{\phi}=(H\cap \fix{\phi})*(K\cap \fix{\phi})$, where $r(K\cap \fix{\phi})=1$, \item there exist non-trivial elements $y\in F$,
$h,h'\in H$, such that $K=\langle y\rangle$, $y\phi =h'y$, $h$ is not a proper power, $\fix{\phi}=(H\cap \fix{\phi})*\langle y^{-1}hy
\rangle$ and $h\phi =h'hh'^{-1}$.
\end{apartats}
\end{Thm}

In the maximal rank case we can compare with Theorem~\ref{maxrank}. In this situation, (i) never happens (by Bestvina-Handel Theorem,
see~\cite{BH}), (ii) corresponds to~\ref{maxrank}(i) (since one can easily deduce that $K$ must be cyclic), and (iii) corresponds
to~\ref{maxrank}(ii) (with the extra information that $h$ can be chosen to be fixed, and then $h'$ to be a power of $h$).

Theorem~\ref{mainconnex} provides an local description of 1-auto-fixed subgroups. Using induction, it is easy to deduce
Theorem~\ref{cormain}, our main result, providing a global description which may be of greater interest to the more general reader.

\begin{Thm}\label{cormain}
Let $F$ be a non-trivial finitely generated free group and let $\phi\in Aut(F)$ with $\fix{\phi} \neq 1$. Then, there exist integers
$r,s\geq 0$, $\phi$-invariant non-trivial subgroups $K_1, \ldots ,K_r \leq F$, primitive elements $y_1, \ldots ,y_s \in F$, a subgroup
$L\leq F$, and elements $1\neq h'_j \in H_j =K_1*\cdots *K_r*\langle y_1,\ldots ,y_j \rangle$, $j=0,\ldots ,s-1$, such that
 $$
F=K_1*\cdots *K_r*\langle y_1,\ldots ,y_s \rangle *L
 $$
and $y_j \phi =h'_{j-1}y_j$ for $j=1,\ldots ,s$; moreover,
 $$
\fix{\phi}=\langle w_1, \ldots ,w_r, y_1^{-1}h_0 y_1, \ldots ,y_s^{-1}h_{s-1}y_s \rangle
 $$
for some non-proper powers $1\neq w_i\in K_i$ and some $1\neq h_j\in H_j$ such that $h_j\phi =h_j' h_j h_j'^{-1}$, $i=1,\ldots ,r$,
$j=0,\ldots ,s-1$.
\end{Thm}

The paper is structured as follows.

In section~\ref{prel} we prove some easy but useful lemmas, that will be needed later. Then, Theorems~\ref{mainconnex}
and~\ref{cormain} will be illustrated with some explicit examples, pointing out some interesting aspects of the results.

The main tool we use is the theory of relative train track maps, developed by Bestvina and Handel~\cite{BH} and extended later by
Bestvina, Feighn and Handel~\cite{BFH}. Due to the technical nature of the paper, we try to provide as self contained a proof as
possible, by writing the preliminary sections~\ref{s BH}, \ref{freefac} and~\ref{improved rtt}, providing a complete introduction and
recalling the parts of the theory we use, adapted to our needs.

Concretely, we make use of the alternative but equivalent formulation of relative train track maps in terms of free groupoids and
their morphisms, due to Dicks and Ventura~\cite{DV}. In section~\ref{s BH} we recall the Bestvina Handel Theory of relative train
track maps as well as we establish the terminology for free groupoids and their morphisms. For the reader unfamiliar with free
groupoids and groupoid morphisms we note that they are the objects and morphisms obtained from applying the $\pi_1$ functor to the
category of graphs and continuous maps relative to a set of vertices rather than a single vertex. Thus, the language of graphs (eg.
connectedness, components, paths, etc.) carries over to the context of free groupoids with the difference that homotopies relative to
the vertex set are equalities in the groupoid. For more details see~\cite{DV} and~\cite{H}.

Section~\ref{freefac} reviews free factor systems and is taken mostly from~\cite{BFH}, while in section~\ref{improved rtt} we recall
the improved relative train tracks from~\cite{BFH} and we prove Theorem~\ref{cbfh}, one of the crucial pieces of the main argument.

The main argument of the paper is the proof of Theorem~\ref{main} in section~\ref{s main theorem}. This result provides a simultaneous
description of all the eigengroups for a given automorphism of $F$, and of its fixed subgroup among them. The proof is a long and
technical inductive argument with several cases, all in the context of graphs, groupoids and relative train track maps. Along it we
need to deal with outer automorphisms instead of automorphisms, and with free factor systems instead of free factors. Finally, in
section~\ref{applications}, Theorem~\ref{mainconnex} as well as our main result Theorem~\ref{cormain}, will be easily deduced from
Theorem~\ref{main}.

From the formal point of view, Theorem~\ref{main} is stronger than Theorem~\ref{cormain}. However, we present the last one as the main
result in the paper because of its greater simplicity and algebraic usefulness, and because of the greater technicality involved in
Theorem~\ref{main}.

One observation that we would like to make is that although our results are stated and proved for automorphisms, we strongly expect
that Theorems~\ref{main}, \ref{mainconnex} and \ref{cormain} hold for monomorphisms of free groups. In fact, as \cite{DV} applies to
monomorphisms, the only place we essentially use the fact that we are dealing with automorphisms is in section~\ref{improved rtt},
Theorems~\ref{bfh} and~\ref{cbfh}, which in turn rely on \cite{BFH}. We are confident that \cite{BFH} holds for monomorphisms (with
minor technical changes) using essentially the same proofs. However, the work needed to verify this seems disproportionate to the
potential gain.

As a direct consequence of the results in this work, the same authors have the subsequent paper~\cite{MVPRE}, where
Theorem~\ref{cormain} is used to show that the collection of 1-endo-fixed subgroups is strictly larger than the collection of
1-auto-fixed subgroups if the underlying free group has rank at least 3.

\section{Lemmas and examples}\label{prel}

In this section we prove some simple but useful lemmas. Then, we construct example~\ref{ex1}, which typically illustrates
Theorems~\ref{mainconnex} and~\ref{cormain}, and examples~\ref{ex2} and~\ref{ex3} pointing out some interesting aspects of these two
theorems.

\begin{Lem}\label{iff}
Let $F$ be a finitely generated free group and let $\phi\in Aut(F)$. If $\fix{\phi}$ is contained in a proper free factor $H$ of $F$
then there exists a (proper) $\phi$-invariant free factor $K$ of $F$ such that $\fix{\phi}\leq K\leq H$.
\end{Lem}

\demo Let $K$ be a free factor of $H$ (and so, of $F$) containing $\fix{\phi}$ and having the smallest possible rank. We only have to
show that $K$ is $\phi$-invariant. Clearly, $K\phi$ is also a free factor of $F$ and so, $K\cap K\phi$ is a free factor of $K$. From
the minimality of $r(K)$ and the fact $\fix{\phi}\leq K\cap K\phi$, we deduce $K\cap K\phi=K$. Thus, $K$ is a free factor of $K\phi$.
But these two subgroups have the same rank, so $K=K\phi$. \qed

Consequently, if the fixed subgroup of an automorphism $\phi\in Aut(F)$ is contained in a proper free factor of $F$ then the action of
$\phi$ which gives rise to fixed elements takes place in a proper $\phi$-invariant free factor of $F$ (and this corresponds to
construction~\ref{construct}~(i)).

\begin{Lem}\label{imagey}
Let $F$ be a free group and let $\phi\in Aut(F)$. Let $H$ be a $\phi$-invariant subgroup of $F$. Suppose that $y\in F$ and $h\in H$
are non-trivial elements such that $F=H*\langle y\rangle$ and $y^{-1}hy\in \fix{\phi}$. Then, there exists $h'\in H$ such that $y\phi
=h'y$ and $h\phi =h'hh'^{-1}$. Furthermore, $h'=1$ if, and only if, $y\in \fix{\phi}$.
\end{Lem}

\demo Since $\phi$ is an automorphism, we have
 $$
F=H*\langle y\rangle =H\phi*\langle y\phi\rangle =H*\langle y\phi\rangle.
 $$
Hence, $y\phi =h' y^{\epsilon}h''$ for some $h',h'' \in H$ and $\epsilon =\pm 1$. Now the equation
 $$
y^{-1}hy =(y^{-1}hy)\phi =h''^{-1}y^{-\epsilon}h'^{-1}(h\phi)h' y^{\epsilon}h''
 $$
forces $h''=1$, $\epsilon =1$ and $h=h'^{-1}(h\phi)h'$. Thus, $y\phi =h'y$ and $h\phi =h'hh'^{-1}$. Clearly, $h'=1$ if, and only if,
$y$ is fixed by $\phi$. \qed

It is well known that every non-trivial free factor $H$ of $F$ is \emph{malnormal}, i.e. $H^x \cap H\neq 1$ only when $x\in H$. This
can be reformulated in the following useful way.

\begin{Lem}\label{malnormal}
Let $F$ be a free group, $H$ a free factor of $F$ and $\phi,\varphi\in Aut(F)$ be such that $\phi =\varphi \gamma_x$ for some $x\in
F$. If $H$ is $\phi$-invariant and $H\cap \fix{\varphi}\neq 1$ then $x\in H$, and $H$ is also $\varphi$-invariant.
\end{Lem}

In general, computing fixed subgroups of automorphisms of free groups is a difficult task (although, recently, some authors have
announced algorithms for this purpose, see~\cite{M} or~\cite{L}). It is even difficult, in general, to show that a given subgroup is
really the fixed subgroup of a given automorphism. In the following example, we compute the fixed subgroup of a given automorphism,
and then we see how it fits with the descriptions given by Theorems~\ref{mainconnex} and~\ref{cormain}.

\begin{Emp}\label{ex1}
Let $F$ be the free group of rank 6 freely generated by the letters $\{a,b,c,d,e,f\}$, and consider the automorphism $\phi$ of $F$
given by
 $$
\begin{array}{rclccrcl}
a & \to & a & \qquad & b &\to & ab \\ c & \to & dc & & d & \to & dcd \\ e & \to & (b^{-1}ab[c,d])^t e & & f & \to & bfb,
\end{array}
 $$
where $t$ is an integer. Choosing $t$ appropriately, we claim that
 $$
\fix{\phi}=\langle a, b^{-1}ab, [c,d], e^{-1}b^{-1}ab[c,d]e \rangle.
 $$
Clearly, both $\langle a,b \rangle$ and $\langle c,d \rangle$ are $\phi$-invariant. Automorphisms of the free group of rank 2 are well
understood and it is easy to check that
 $$
\fix{\phi_{\langle a,b,c,d \rangle}}=\fix{\phi_{\langle a,b \rangle}}*\fix{\phi_{\langle c,d \rangle}}=\langle a, b^{-1}ab, [c,d]
\rangle.
 $$
Since $b^{-1}ab[c,d]$ is fixed by $\phi$, the proposition in the Appendix tells us that, for all but finitely many choices of $t$,
 $$
\fix{\phi_{\langle a,b,c,d,e \rangle}}=\langle a, b^{-1}ab, [c,d], e^{-1}b^{-1}ab[c,d]e \rangle.
 $$
We choose and fix one such value for $t$ (in fact, $t=1$ would suffice). It remains to see that $f$ contributes nothing to the fixed
subgroup of $\phi$. To see this, suppose that $w=w\phi$ is a fixed word not contained in $\langle a,b,c,d,e \rangle$ and find a
contradiction. The reduced form of $w$ looks like
 $$
w=u_0 f^{\epsilon_1} \cdots f^{\epsilon_i} u_i f^{\epsilon_{i+1}} \cdots f^{\epsilon_k} u_k,
 $$
where $k\geq 1$, $\epsilon_i=\pm 1$, and the $u_i \in \langle a,b,c,d,e \rangle$ are possibly trivial but otherwise reduced as
written. We then get that
 $$
w\phi=(u_0\phi) b^{\epsilon_1} f^{\epsilon_1} b^{\epsilon_1} \cdots b^{\epsilon_i} f^{\epsilon_i} b^{\epsilon_i} (u_i\phi)
b^{\epsilon_{i+1}} f^{\epsilon_{i+1}} b^{\epsilon_{i+1}} \cdots b^{\epsilon_k }f^{\epsilon_k} b^{\epsilon_k} (u_k\phi).
 $$
Since $w=w\phi$, no pair of $f$'s can cancel in the last expression and, in particular, $(u_0\phi) b^{\epsilon_1}= u_0 $.
Abelianising, we obtain $u_0^{\rm ab}\in \mathbb{Z}^5$ and $\phi_{\langle a,b,c,d,e\rangle}^{\rm \, ab}\in GL_5 (\mathbb{Z})$ such
that $u_0^{\rm ab}\phi_{\langle a,b,c,d,e\rangle}^{\rm \, ab}=u_0^{\rm ab}-\epsilon_1 b^{\rm ab}$. But, from the definition of $\phi$,
it is easy to check that $b^{\rm ab}$ is not in the image of $\phi_{\langle a,b,c,d,e\rangle}^{\rm \, ab}-Id$. This contradiction
proves the claim:
 $$
\fix{\phi}=\fix{\phi_{\langle a,b,c,d,e\rangle}}=\langle a, b^{-1}ab, [c,d],e^{-1}b^{-1}ab[c,d]e \rangle.
 $$
This 1-auto-fixed subgroup of $F$ fits the description of Theorem~\ref{mainconnex} in the following way. The automorphism $\phi$ is in
case (i) with $\langle a,b,c,d,e \rangle$ invariant and containing the whole fixed subgroup of $\phi$. Now, $\phi_{\langle a,b,c,d,e
\rangle}$ is in case (iii) with $\langle a,b,c,d \rangle$ invariant, $y=e$, $h'=(b^{-1}ab[c,d])^t$ and $h=b^{-1}ab[c,d]$. Then,
$\phi_{\langle a,b,c,d \rangle}$ is in case (ii) with the two invariant free factors being $\langle a,b \rangle$ and $\langle c,d
\rangle$. Finally, $\phi_{\langle c,d\rangle}$ has cyclic fixed subgroup, and $\phi_{\langle a,b\rangle}$ is again in case (iii) with
$\langle a\rangle$ invariant, $y=b$ and $h=h'=a$. Note that the last two steps can also be interchanged. This shows how $\fix{\phi}$
can be built up using the basic constructions described in~\ref{construct}.

Also, $\fix{\phi}$ fits the description of Theorem~\ref{cormain} as follows. Take $r=s=2$, $K_1=\langle a \rangle$, $K_2=\langle c,d
\rangle$, $y_1=b$, $y_2=e$ and $L=\langle f \rangle$. Then, taking $w_1=a\in K_1$, $w_2=[c,d]\in K_2$, $h'_0=h_0=a\in K_1*K_2$, and
$h'_1=(b^{-1}ab[c,d])^t$ and $h_1=b^{-1}ab[c,d]$ both in $K_1*K_2*\langle b\rangle$, we have
 $$
\fix{\phi}=\langle w_1, w_2, y_1^{-1}h_0y_1, y_2^{-1}h_1y_2 \rangle =\langle a, [c,d], b^{-1}ab, e^{-1}b^{-1}ab[c,d]e\rangle.
 $$
\end{Emp}

One word of caution is necessary here since the $h$-elements in the example above,
 corresponding to $y_1=b$ and $y_2=e$ (i.e.
$h_0=a$ and $h_1=b^{-1}ab[c,d]$, respectively) are fixed by $\phi$, whereas in Theorem~\ref{cormain} they have a more complicated
image. We might hope, with the correct choice of basis, to always ensure that all the $h$-elements in Theorem~\ref{cormain} are fixed,
as happens in the maximal rank case. This will not be possible in general as the next example shows.

\begin{Emp}\label{ex2}
Let $F$ be the free group of rank 3 freely generated by $\{a,b,c\}$. Consider the elements $g=[a,b]$ and $h=a^2b^2$ of $\langle a,b
\rangle$, which are not proper powers, and let $\phi$ be the automorphism of $F$ given by $a\mapsto g^{-1}ag$, $b\mapsto g^{-1}bg$,
$c\mapsto g^{-1}hc$. Using cancellation arguments, it is not difficult to see that $$\fix{\phi}=\langle g, c^{-1}hc \rangle.$$ We
claim that $F$ does not have a basis for $\phi$ as in Theorem~\ref{cormain} with all the $h$-elements being fixed.

Firstly, by looking at the abelianisation of $F$, it is clear that the subgroup $\fix{\phi}=\langle g, c^{-1}hc \rangle$ contains no
primitive elements  of $F$. So, if we obtain a basis for $F$ and a description of $\fix{\phi}$ as in Theorem~\ref{cormain}, the
subgroup $K_1$ must have rank $2$ and hence, $L=1$ and $\fix{\phi}=\langle w_1, {y_1}^{-1} h_0 y_1 \rangle$, where $w_1, h_0 \in K_1$
are not proper powers. Suppose now that $h_0$ is fixed and find a contradiction.

In this situation, $h_0={w_1}^{\pm 1}$. Then, the normal subgroup generated by $w_1$ coincides with the normal subgroup generated by
$g$ and $h$. Consider the one relator group $G=\langle a,b,c \ | \ w_1 \rangle =\langle a,b,c \ |  \ g,h \rangle$. Since $w_1$ is not
a proper power, $G$ is torsion free (see Proposition~5.18 in~\cite{LS}). However, looking at the second presentation, it is clear that
$ab$ is an order two element of $G$. This contradiction demonstrates that $\fix{\phi}$ cannot be described as in Theorem~\ref{cormain}
with fixed $h$-elements.
\end{Emp}

Finally, we remark that Theorem~\ref{cormain} provides a description of 1-auto-fixed subgroups of $F$, but it does not give a
characterisation of those subgroups. The following is an example of a subgroup which agrees with the description given, while it is
not 1-auto-fixed.

\begin{Emp}\label{ex3}
Let $F$ be the free group of rank 3 freely generated by $\{a,b,c\}$ and consider the subgroup $H=\langle a, b^{-1}ab, c^{-1}bc\rangle
\leq F$. Clearly, $H$ agrees with the description given in Theorem~\ref{cormain}. Suppose that $H=\fix{\phi}$ for some $\phi \in
Aut(F)$. It is straightforward to verify that $\langle a,b\rangle $ is the smallest free factor of $F$ containing the subgroup
$\langle a,b^{-1}ab\rangle$. Hence, $\langle a,b\rangle$ is $\phi$-invariant, that is, $\langle a,b\rangle =\langle a, b\phi\rangle$.
Thus, $b\phi =a^r b^{\epsilon}a^s$ where $\epsilon =\pm 1$ and $r,s$ are integers. But the elements $a$ and $b^{-1}ab$ being fixed
imply $a\phi =a$ and $b\phi =a^r b$ for some $r\neq 0$. Now $F=\langle a,b\rangle *\langle c\rangle $ where $\langle a,b\rangle $ is
$\phi$-invariant, $c\not\in \fix{\phi}$ and $c^{-1} b c \in \fix{\phi}$. Thus we may apply Lemma~\ref{imagey} to conclude that
$a^rb=b\phi =h'bh'^{-1}$ for some $1\neq h'\in \langle a,b\rangle $. This contradiction shows that $H$ is not 1-auto-fixed. In fact,
it is easy to see that any endomorphism of $F$ fixing $H$ has to actually be an automorphism, and has to fix $\langle a, b,
c^{-1}bc\rangle$ too. Hence, $H$ is not even the intersection of an arbitrary family of 1-endo-fixed subgroups of $F$.
\end{Emp}

\section{Bestvina-Handel theory}\label{s BH}

In this paper, we will mostly work with finite graphs viewed as combinatorial objects. We will make use of the standard concepts and
terminology for graphs, including the notions of rank and reduced rank of a graph, connectedness, subgraph and the core of a graph,
graph morphism, (formal) path, trivial path, reduced (or normal) form of a path, etc. We refer the reader to~I.1.2 and~I.1.8
in~\cite{DV} for the precise definitions. For example, if $Z$ is a graph, $VZ$ denotes the vertex set, $EZ$ the edge set, $\iota$ and
$\tau$ are the incidence functions, $\cdot$ denotes the concatenation between paths (when defined), etc. If $Z$ is connected (and only
then), $r(Z)$ is the rank of $Z$, the rank of the fundamental group of $Z$. The reduced rank of $Z$, $\tilde{r}(Z)$, will be the
reduced rank of its fundamental group; the maximum of $r(Z)-1$ and zero. If $Z$ is not connected $r(Z)$ is not defined, and
$\tilde{r}(Z)$ is defined to be the sum of the reduced ranks of its components.

The excellent paper~\cite{BH} deals with graphs in the topological point of view, and~\cite{DV} contains a complete reformulation of
it into the language of groupoids (and, additionally, many details previously left to the reader are meticulously verified). We found
the language used in~\cite{DV} more convenient to our purposes, so we mostly will refer to it. However, we will also emphasize the
corresponding concepts in~\cite{BH} in order to make the arguments more clear.

Let $Z$ be a graph. In~\cite{BH} homotopy equivalences of graphs, $\beta \colon Z\rightarrow Z$, are considered. For the purposes of
that work (and the present one), all the relevant data of such a homotopy equivalence is what it induces at the fundamental group
level. So, it is useful to forget the topological structure of the graph $Z$ and to just think of it as a combinatorial object. In
this setting, $\beta$ must be thought of as a formal map sending edges (and so paths) to paths, and respecting the incidences in $Z$.

A good language to formalise this is the language of groupoids.

\begin{Dfn}\label{defs1}
A \emph{groupoid} is a small category in which the objects, also called \emph{vertices}, are identified with the identity morphisms,
and every morphism is invertible (see~\cite{DV} or~\cite{H} for more details). For example, every group is a groupoid with a single
vertex (and so, the operation is totally defined). Note that every groupoid has also the structure of a graph.

It is easy to check that if $Z$ is a graph, then the set of paths in $Z$ modulo reduction, together with the natural structure coming
from $Z$, form a groupoid. It is called the \emph{fundamental groupoid of} $Z$ and denoted $\pi Z$. To simplify the notation, we
identify every class of equivalent paths with the unique reduced path it contains, called its \emph{normal form}. For every $u, v\in
VZ$, $\pi Z(u,v)$ denotes the set of paths in $\pi Z$ from $u$ to $v$. In particular, $\pi Z(u,u)$, which is a free subgroup of $\pi
Z$, is the (combinatorial) \emph{fundamental group} of $Z$ at $u$. A path $p\in \pi Z(u,v)$ defines an isomorphism, $\gamma_p \colon
\pi Z(u,u) \mapsto \pi Z(v,v)$ by sending $x\in \pi Z(u,u)$ to $p^{-1}xp \in \pi Z(v,v)$. Clearly, the inverse path, $p^{-1}$, defines
the inverse isomorphism $\gamma_{p^{-1}}\colon \pi Z(v,v) \mapsto \pi Z(u,u)$. Thus, in the case where $Z$ is connected, all the
groups $\pi Z(v,v)$ are isomorphic (in fact, all of them are free of rank $r(Z)$).

A graph $T$ is a \emph{tree} if, and only if, for every $u,v\in VT$, $\pi T(u,v)$ contains a unique element, called the
\emph{geodesic} from $u$ to $v$ and denoted $T[u,v]$.

We say that a groupoid $G$ is \emph{free} when it is of the form $\pi Z$ for some subgraph $Z$ of $G$. In this case, we also say that
$Z$ is a \emph{basis} of $G$. It can be proved that the Nielsen-Schreier theorem is also valid for groupoids, (see Theorem~I.3.10
in~\cite{DV}). So, every subgroupoid $H$ of $G=\pi Z$ is free and has a basis $B$, which is a graph whose vertices are vertices of
$Z$, and whose edges are elements of $B\leq \pi Z$. In general, $B$ can be disconnected, even when $Z$ is connected. Furthermore, if
$B$ and $Z$ are connected, $r(B)$ can be bigger than $r(Z)$, as in the case of free groups.
\end{Dfn}

\begin{Dfn}\label{defs2}
In the language of groupoids, continuous maps $\beta\colon Z\rightarrow Z$ become groupoid morphisms $\pi \beta \colon \pi Z
\rightarrow \pi Z$, also referred to as \emph{self-maps} of $Z$, simply denoted $\beta \colon Z\rightarrow Z$. And homotopy
equivalences become \emph{equivalences of graphs} that is, self-maps $\beta$ such that, for every $v\in VZ$, the restriction
$\beta_v\colon \pi Z(v,v)\to \pi Z(v\beta,v\beta)$ is an isomorphism of groups.

A path $p\in \pi Z$ is called \emph{$\beta$-fixed} when $p\beta =p$. We include here the case where $p$ is trivial, i.e. a path
consisting of a single vertex $v$, in which case we say that $v$ is $\beta$-fixed. Clearly, the set of $\beta$-fixed paths,
$\fix{\beta}=\{ p\in \pi Z : p\beta =p\}$, is a subgroupoid of $\pi Z$.
\end{Dfn}

In order to study fixed subgroups of automorphisms of free groups, what Bestvina and Handel really did in~\cite{BH} was to study the
fixed subgroupoid of a given equivalence of graphs. In their Proposition~6.3, they constructed a basis for $\fix{\beta}$, with
technical assumptions on $\beta$.

\begin{Dfn}\label{defs3}
Recall that, throughout the paper, $F$ denotes a finitely generated free group.

Following~\cite{BH}, a \emph{marked graph} is a triple $(Z,v_0,\eta)$, where $Z$ is a connected graph, $v_0\in VZ$ is a vertex called
\emph{basepoint}, and $\eta \colon F \rightarrow \pi Z(v_0,v_0)$ is a group isomorphism, called a \emph{marking}.

Let $\beta \colon Z\rightarrow Z$ be an equivalence of a marked graph $(Z, v_0, \eta)$. For every path $p\in \pi Z(v_0\beta ,v_0)$,
one can consider the \emph{induced} automorphism $\phi_{\beta, p,\eta}$ of $F$ defined by $x\mapsto (p^{-1}\cdot x\eta \beta \cdot p
)\eta^{-1}$ (see the commuting diagram below). It will be denoted $\phi_{\beta,p}$ whenever the marking is understood from the
context.

 $$
\xymatrix{ F \ar[dddd]_{\phi_{\beta, p,\eta}} \ar[rr]^{\eta \qquad } & & \pi Z(v_0, v_0) \ar[dd]^{\beta} \\ \\ &  &\pi Z(v_0\beta ,
v_0 \beta) \ar[dd]^{\gamma_p} \\ \\ F \ar[rr]^{\eta \qquad} & & \pi Z(v_0, v_0)\\ }
 $$

In~\cite{BH} it is noted that, changing the path $p$, the automorphism $\phi_{\beta, p}$ is composed with an inner automorphism.
Conversely, for every $y\in F$, it is clear that $p\cdot y\eta\in \pi Z(v_0\beta, v_0)$ and $\phi_{\beta, p\cdot y\eta}=\phi_{\beta,
p}\gamma_y$. Hence, $\{ \phi_{\beta, p} : p\in \pi Z(v_0 \beta, v_0) \}$ is an outer automorphism of $F$ determined by $\beta$ and
denoted $\Phi_\beta$. It is said that $\beta$ \emph{induces} or \emph{represents} the outer automorphism $\Phi_{\beta}$ with respect
to the given marking, $\eta$. However, we shall suppress $\eta$ from the notation if there is no risk of confusion.
\end{Dfn}

In the particularly simple case where $v_0$ is $\beta$-fixed and $p=v_0$ is a trivial path, we have $\phi_{\beta, v_0}=\eta
\beta_{v_0} \eta^{-1}$ and, identifying $F$ with $\pi Z(v_0,v_0)$ via $\eta$, we can abuse notation and write $\phi_{\beta,
v_0}=\beta_{v_0}$. The Fixed-point Lemma below ensures that, when the fixed subgroup is non-cyclic, such a simple marking can always
be chosen.

\begin{Lem}[Fixed-point lemma]\label{fpl}
Let $\beta\colon Z\rightarrow Z$ be an equivalence of a mark\-ed graph $(Z,v_0,\eta)$, such that all $\beta$-fixed points are vertices
(see~I.5.1 in~\cite{DV}). For every path $p\in \pi Z(v_0\beta, v_0)$ with $r(\fix{\phi_{\beta,p}})\geq 2$ there exist a $\beta$-fixed
vertex $v\in VZ$ and a path $q\in \pi Z(v, v_0)$ such that $q\beta =q\cdot p^{-1}$. Moreover,
\begin{apartats}
\item $\phi_{\beta,v,\eta'}=\phi_{\beta ,p, \eta}$, where
$\eta'\colon F\to \pi Z(v,v)$, $x\mapsto q\cdot x\eta \cdot q^{-1}$, \item $x\phi_{\beta,p}=(q^{-1}\cdot q\beta\cdot x\eta\beta \cdot
(q\beta)^{-1}\cdot q)\eta^{-1}$ for every $x\in F$, \item $\fix{\phi_{\beta,p}}=(q^{-1}\cdot \fix{\beta_v}\cdot q)\eta^{-1}$.
\end{apartats}
\end{Lem}

\demo Corollary~2.2 in~\cite{BH} or Lemma~I.5.4 in~\cite{DV} give us the required vertex $v$ and path $q$. By looking at the following
commuting diagram, we see (i). Equalities (ii) and (iii) are straightforward to verify. \qed

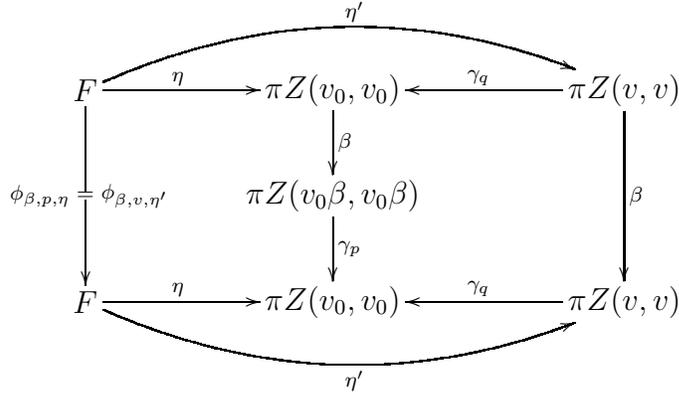
 $$
\xymatrix{ F \ar[dd]_{\phi_{\beta, p,\eta}\ }|=^{{ \ \phi_{\beta, v,\eta'}}} \ar[rr]^{\eta \qquad } \ar@/^2pc/[rrrr]^{\eta'} & & \pi
Z(v_0, v_0) \ar[d]^{\beta}  & &\pi Z(v,v) \ar[ll]_{\gamma_{q}} \ar[dd]^{\beta}\\ &  & \pi Z(v_0\beta , v_0 \beta) \ar[d]^{\gamma_p} \\
F \ar[rr]^{\eta \qquad} \ar@/_2pc/[rrrr]_{\eta'} & & \pi Z(v_0, v_0) & & \pi Z(v,v) \ar[ll]_{\gamma_{q}} \\ }
 $$

\begin{Dfn}\label{defs4}
Let $\beta$ be an equivalence of a finite, connected, marked, core graph $(Z,v_0,\eta)$, and let $Z_0$ be a maximal proper
$\beta$-invariant subgraph of $Z$ (note that $Z_0$ is not necessarily connected, and if $EZ\neq \emptyset$ then $VZ_0 =VZ$). As with
$\beta$, we say that the triple $(\beta,Z,Z_0)$ \emph{represents} $\Phi_{\beta}$. In fact, in~\cite{BH}, the authors constructed a
full filtration of $Z$ by proper $\beta$-invariant subgraphs, and $(\beta,Z,Z_0)$ corresponds to the top stratum. One of the main
simplifications introduced in~\cite{DV} is to work only with this top stratum (that is, with a maximal proper $\beta$-invariant
subgraph) instead of working with the full filtration, and then do the arguments inductively (this idea originally came from Gaboriau,
Levitt and Lustig,~\cite{GLL}).

In~\cite{BH}, an irreducible matrix is associated to $(\beta,Z,Z_0)$. This matrix, denoted $[\beta/Z_0]$ in~\cite{DV}, plays a central
role in the arguments. By Perron-Frobenius Theorem (see Theorem~1.5 in~\cite{BH} or section~II.1 in~\cite{DV}), the spectral radius of
$[\beta/Z_0]$ is 0, 1 or $>1$; we will refer to these three possibilities by saying that $(\beta,Z,Z_0)$ is \emph{null}, is
\emph{level} or is \emph{exponential}, respectively. This spectral radius is denoted $\lambda$ in~\cite{BH}, and $PF(\beta/Z_0)$
in~\cite{DV}.
\end{Dfn}

One of the main results in~\cite{BH} states that every  $\Phi \in Out(F)$ can be represented by some $(\beta,Z,Z_0)$ with several
extra good properties; these are the \emph{stable relative train tracks}. The corresponding notion in~\cite{DV} is that of
\emph{minimal representatives}. Concretely, in section~IV.1 of~\cite{DV} it is proved that, given an outer automorphism $\Phi \in
Out(F)$, there exist a finite, connected, marked, core graph $(Z,v_0,\eta)$, a self-map $\beta\colon Z\rightarrow Z$ representing
$\Phi$, and a maximal, proper, $\beta$-invariant subgraph $Z_0$ of $Z$ such that $(\beta, Z, Z_0)$ is a minimal representative of
$\Phi$ (in fact, an analogous result is proven for injective endomorphisms, but here we are only interested in the bijective case).
The precise definition of a minimal representative (or a relative train track) is not important here. All we will use about them is
the existence of such representatives and their main properties. A first consequence of the definition is that, for every minimal
representative $(\beta,Z,Z_0)$ of $\Phi$, the value $PF(\beta/Z_0)$ is the minimum possible among those of all the representatives of
$\Phi$. The other useful properties of minimal representatives are expressed in Theorem~IV.5.1. We rewrite this statement here in a
slightly different form, distinguishing the three possibilities for the stratum (see the corresponding proof in~\cite{DV}), and
restricting the attention to isomorphisms (for which $(\beta,Z,Z_0)$ cannot be null).

\begin{Thm}[IV.5.1,~\cite{DV}]\label{threpr}
Let\, $F$ be\, a\, non-trivial\, finitely\, generated\, free\, group, and let $\Phi \in Out(F)$. Then, there exists an equivalence
$\beta$ of a non-empty, finite, connected, marked, core graph $(Z,v_0,\eta)$, such that $\Phi_{\beta}=\Phi$ and
\begin{apartats}
\item all $\beta$-fixed points are vertices of $Z$ (see
Definition~I.5.1 in~\cite{DV}),
\item $Z$ has a proper $\beta$-invariant subgraph $Z_0$ on which
$\beta$ induces a self-map $\beta_0$,
\item $\fix{\beta}$ has a basis $B$ containing (as a subgraph) a
basis $B_0$ of $\fix{\beta_0}$,
\item $(\beta, Z,Z_0)$ is not null,
\item if $(\beta, Z,Z_0)$ is level then either $B=B_0$ or
$B\setminus B_0 =\{ e \}$, where $e$ is a $\beta$-fixed edge such that $EZ\setminus EZ_0=\{ e\}$,
\item if $(\beta,Z,Z_0)$ is exponential then either $B=B_0$ or
$B\setminus B_0 =\{ p \}$ where $p\in \pi Z$ is such that either some edge $e\in EZ\setminus EZ_0$ $Z$-occurs in $p$ only once, or
every edge $e\in EZ\setminus EZ_0$ $Z$-occurs in $p$ exactly twice.
\end{apartats}
In fact, every minimal representative $(\beta, Z,Z_0)$ of $\Phi$ satisfies properties (i) to (vi).
\end{Thm}

We want also to remark on the following two technical differences between~\cite{BH} and~\cite{DV}. The first is that in~\cite{BH} the
graph $Z_0$ does not contain isolated vertices (it is defined as the closure of certain set of edges) while in~\cite{DV} it contains
all the vertices of $Z$ (it is defined as a maximal but not necessarily connected, invariant subgraph); up to isolated vertices, they
are the same subgraph of $Z$. The second technical difference is about fixed paths: in~\cite{BH} they consider Nielsen (i.e. fixed)
paths possibly crossing partial edges, while this has no sense in~\cite{DV}; but, using subdivisions at non-vertex fixed points
(see~III.2 in~\cite{DV}), this situation is modelled in~\cite{DV} by assuming that ``all fixed points are vertices", and then looking
only at fixed paths beginning and ending on vertices.

Finally, let us discuss the main result in~\cite{BH} and write it in a slightly different form. With the notation of
Theorem~\ref{threpr}, Bestvina and Handel proved that $\tilde{r}(B)\leq \tilde{r}(F)$ and, as an immediate consequence, they obtained
their main result (Theorem~6.1 in~\cite{BH}), saying that $r(\fix{\phi})\leq r(F)$ for every $\phi\in Aut(F)$. By analyzing the
algebraic meaning of the components of $B$, one can obtain a non-connected version of Bestvina-Handel Theorem, not explicitly stated
in~\cite{BH} (see Theorem~\ref{bh no connex} below or take $H=F$ in Theorem~IV.5.5 of~\cite{DV}). As noted in~\cite{GLL}, this version
of the result is only superficially stronger. However, we will dedicate the rest of the present section to state this non-connected
version, for later use.

For the rest of the section, let $\Phi \in Out(F)$, and choose $\phi\in \Phi$.

Consider the following two equivalence relations in the sets $F$ and $\Phi$, respectively. We say that $y_1, y_2\in F$ are {\it
Reidemeister equivalent}, denoted $y_1\sim_{\phi} y_2$, if $y_2=(c\phi)^{-1}y_1 c$ for some $c\in F$. And we say that $\phi
\gamma_{y_1}, \phi \gamma_{y_2}\in \Phi$ are {\it isogredient}, denoted $\phi \gamma_{y_1} \sim \phi \gamma_{y_2}$, if they are equal
up to conjugation (in $Aut(F)$) by an inner automorphism. That is, when $\phi \gamma_{y_2} =\gamma_{c}^{-1}\phi \gamma_{y_1}
\gamma_c=\gamma_{c^{-1}}\phi \gamma_{y_1 c}$ for some $\gamma_c\in Inn(F)$. (We prefer to use the word ``isogredience" instead of, for
example, ``similarity", used by other authors; the reason is historical since, although in a different context, Nielsen already dealt
with the same concept under the name ``isogredience").

Note that the isogredience relation does not depend on the chosen $\phi$, while Reidemeister relation does, like the particular
bijection $y\mapsto \phi \gamma_y$ from $F$ to $\Phi$. The following lemma clarifies the relationship between these two equivalence
relations, and the eigengroups of $\phi$.

\begin{Lem}\label{nc1}
Let $F$ be a non-cyclic finitely generated free group and $y_1,y_2 \in F$. Let $\Phi\in Out(F)$ and choose $\phi\in \Phi$ such that
$r(\fix{\phi \gamma_{y_1}})\geq 2$. The following are equivalent:
\begin{apartats}
\item[\upshape{(a)}] $y_1 \sim_{\phi} y_2$, \item[\upshape{(b)}]
$\phi \gamma_{y_1} \sim \phi \gamma_{y_2}$, \item[\upshape{(c)}] $[[\fix{\phi \gamma_{y_1}}]]=[[\fix{\phi \gamma_{y_2}}]]$.
\end{apartats}
\end{Lem}

\demo (a)~$\Rightarrow$~(b). By (a), there exists $c\in F$ such that $y_2 =(c\phi)^{-1} y_1 c$. The same element $c$ satisfies the
equality $\phi \gamma_{y_2} =\phi \gamma_{(c \phi)^{-1}y_1 c}=\phi \gamma_{c^{-1}\phi}\, \gamma_{y_1 c}=\gamma_{c^{-1}}\phi
\gamma_{y_1 c}$ in $Aut(F)$. So, $\phi \gamma_{y_1} \sim \phi \gamma_{y_2}$.

(b)~$\Rightarrow$~(a). By (b), there exists $c\in F$ such that $\phi \gamma_{y_2}= \gamma_{c}^{-1}\phi \gamma_{y_1}\gamma_c$. That is,
for every $x\in F$,
 $$
y_2^{-1} (x\phi) y_2 =c^{-1}y_1^{-1}((cxc^{-1})\phi )y_1 c =c^{-1} y_1^{-1} (c\phi) (x\phi) (c\phi)^{-1} y_1 c.
 $$
So, $(c\phi)^{-1} y_1 c y_2^{-1}$ commutes with $x\phi$, for every $x\in F$. Since (the $\phi$-image of) $F$ is not cyclic, and
centralizers of non-trivial elements in a free group are always cyclic, we deduce that $(c\phi)^{-1} y_1 c y_2^{-1}=1$. Hence, $y_1
\sim_{\phi} y_2$.

(a)~$\Rightarrow$~(c). By (a), there exists $c\in F$ such that $y_2 =(c\phi)^{-1} y_1 c$. The same element $c$ satisfies the equality
$\fix{\phi \gamma_{y_2}}=(\fix{\phi \gamma_{y_1}})^c$.

(c)~$\Rightarrow$~(a). Suppose that $\fix{\phi \gamma_{y_2}}=(\fix{\phi \gamma_{y_1}})^c$ for some $c\in F$. Then, for every $x\in
\fix{\phi \gamma_{y_1}}$, we know that $x^c\in \fix{\phi \gamma_{y_2}}$ and so $y_1^{-1} (c\phi) y_2 c^{-1}$ commutes with $x$. Since
$r(\fix{\phi \gamma_{y_1}})\geq 2$, we deduce that $y_1^{-1} (c\phi) y_2 c^{-1}=1$ and hence, $y_1\sim_{\phi} y_2$. \qed

Note that the hypothesis $r(\fix{\phi \gamma_{y_1}})\geq 2$ is only used in the proof of the implication (c)~$\Rightarrow$~(a). So,
(a)~$\Leftrightarrow$~(b)~$\Rightarrow$~(c) are valid for an arbitrary $\phi\in \Phi$.

Now, let $\beta \colon Z\to Z$ be an equivalence of a marked graph $(Z,v_0,\eta)$ such that all $\beta$-fixed points are vertices, and
$\Phi_{\beta}=\Phi$. Choose $p\in \pi Z(v_0\beta, v_0)$ such that $\phi_{\beta, p}=\phi$. For every $y\in F$ with $2\leq r(\fix{\phi
\gamma_y})=r(\fix{\phi_{\beta,p\cdot y\eta}})$, an application of the Fixed-point Lemma~\ref{fpl} ensures the existence of a
$\beta$-fixed vertex $v\in VZ$ and a path $q\in \pi Z(v,v_0)$ such that $q\beta =q\cdot (y\eta)^{-1}\cdot p^{-1}$ and
$\fix{\phi_{\beta,p\cdot y\eta}}=(q^{-1}\cdot \fix{\beta_v}\cdot q)\eta^{-1}$. The following lemma provides the graph theoretic
equivalent of the statements in Lemma~\ref{nc1} above.

\begin{Lem}\label{nc2}
Let $F$ be a non-cyclic finitely generated free group and $y_1,y_2 \in F$. Let $\Phi\in Out(F)$ and choose some $\phi\in \Phi$. Let
$\beta \colon Z\to Z$ be an equivalence of a marked graph $(Z,v_0,\eta)$ such that all $\beta$-fixed points are vertices, and
$\Phi_{\beta}=\Phi$; let $p\in \pi Z(v_0 \beta, v_0)$ be such that $\phi_{\beta,p}=\phi$. Suppose that $r(\fix{\phi \gamma_{y_1}})\geq
2$ and also that $r(\fix{\phi \gamma_{y_2}})\geq 2$. Let $v_1\in VZ$, $q_1\in \pi Z(v_1, v_0)$, and $v_2\in VZ$, $q_2\in \pi
Z(v_2,v_0)$ be as in the previous paragraph for $\phi \gamma_{y_1}$ and $\phi \gamma_{y_2}$, respectively. Then, conditions \emph{(a),
(b)} and \emph{(c)} in Lemma~\ref{nc1} are also equivalent to
\begin{apartats}
\item[\upshape{(d)}] there exists a $\beta$-fixed path $q\in \pi
Z(v_1, v_2)$.
\end{apartats}
\end{Lem}

\demo (a)~$\Rightarrow$~(d). By (a), there exists $c\in F$ such that $y_2 =(c\phi)^{-1} y_1 c$. Then, it is straightforward to verify
that the path $q=q_1 \cdot c\eta \cdot q_2^{-1}\in \pi Z(v_1, v_2)$ is $\beta$-fixed.

(d)~$\Rightarrow$~(a). Consider the closed path $q_1^{-1}\cdot q\cdot q_2\in \pi Z(v_0,v_0)$ and the element $c=(q_1^{-1}\cdot q\cdot
q_2)\eta^{-1}\in F$. It is also straightforward to verify that $(c\phi)^{-1}y_1 c=y_2$. Hence, $y_1 \sim_{\phi} y_2$. \qed

These equivalent statements can be expressed in the following way which will be very useful later, in the development of our main
argument.

\begin{Prop}\label{corbij}
Let $F$ be a non-cyclic finitely generated free group, and let $\Phi\in Out(F)$. Let $\beta \colon Z\to Z$ be an equivalence of a
marked graph $(Z,v_0,\eta)$ such that all $\beta$-fixed points are vertices, and $\Phi_{\beta}=\Phi$. Let $B$ be a basis for
$\fix{\beta}$. There is a bijection from the set of isogredience classes in $\Phi$ with non-cyclic fixed subgroup, to the set of
components of $B$ with rank $\geq 2$.
\end{Prop}

\demo We choose and fix a path $p\in \pi Z(v_0 \beta, v_0)$ and let $\phi=\phi_{\beta, p, \eta} \in \Phi$ as in
Definitions~\ref{defs3}.

For every $y\in F$ with $r(\fix{\phi \gamma_y})\geq 2$, the Fixed-point Lemma~\ref{fpl} gives us a $\beta$-fixed vertex $v\in VZ$ and
a path $q\in \pi Z(v,v_0)$ such that the following diagram commutes:
 $$
\xymatrix{ F \ar[dd]_{\phi \gamma_{y} \ }|=^{{\ \phi_{\beta , p \cdot y \eta}}} \ar[rr]^{\eta \qquad } & & \pi Z(v_0, v_0)
\ar[d]^{\beta} & &\pi Z(v,v) \ar[ll]_{\gamma_{q}} \ar[dd]^{\beta} \\ &  & \pi Z(v_0 \beta , v_0 \beta) \ar[d]^{\gamma_{p \cdot y
\eta}} \\ F \ar[rr]^{\eta \qquad}  & & \pi Z(v_0, v_0) & & \pi Z(v,v) \ar[ll]_{\gamma_{q}} \\ }
 $$
Let $B_y$ be the component of $B$ containing $v$. We have
 $$
\fix{\phi \gamma_y}=\fix{\phi_{\beta,p\cdot y\eta}}=(q^{-1}\cdot \fix{\beta_v}\cdot q)\eta^{-1} =(q^{-1}\cdot \pi B_y(v,v)\cdot
q)\eta^{-1}
 $$
and so, $r(B_y)=r(\fix{\phi \gamma_y})\geq 2$. This observation together with the implication (b)~$\Rightarrow$~(d) above show that
the map sending the isogredience class of $\phi \gamma_y$ to $B_y$ is well-defined. The implication (d)~$\Rightarrow$~(b) proves that
this map is injective.

Let $C$ be a component of $B$ with rank at least 2. Take a ($\beta$-fixed) vertex $v\in VC$ and an arbitrary path $q\in \pi Z(v,
v_0)$. Then, $q\beta\in \pi Z(v, v_0\beta)$ and we have $q\beta =q\cdot (q^{-1}\cdot (q\beta)\cdot p)\cdot p^{-1}=q\cdot
(y\eta)^{-1}\cdot p^{-1}$ for suitable $y\in F$. We then get a commutative diagram as above. Hence, $B_y=C$ and $r(\fix{\phi
\gamma_y})=r(C)\geq 2$. So, our map is bijective, when considered from the set of isogredience classes in $\Phi$ with non-cyclic fixed
subgroup, to the set of connected components of $B$ with rank $\geq 2$. \qed

Since groups and graphs of ranks 0 and 1 have reduced rank 0, the bijection in Proposition~\ref{corbij} allows us to express the
Bestvina-Handel Theorem in the following way (see Theorem~6.1 in~\cite{BH}, and Theorem~IV.5.5 in~\cite{DV} with $H=F$):

\begin{Thm}[Bestvina-Handel,~\cite{BH}]\label{bh no connex}
Let $F$ be a finitely generated free group, let $\Phi\in Out(F)$, and let $\phi\in \Phi$. Then,
 $$
\sum_{y\in F/\sim_{\phi}} \tilde{r}(\fix{\phi \gamma_y}) =\sum_{\varphi\in \Phi/\sim} \tilde{r}(\fix{\varphi}) =
\sum_{[[\mbox{\scriptsize\rm Fix } \phi \gamma_y]]} \tilde{r}([[\fix{\phi \gamma_y}]])\leq \tilde{r}(F).
 $$
In particular, $r(\fix{\phi})\leq r(F)$, and $\phi$ has at most $\tilde{r}(F)$ non-cyclic conjugacy classes of eigengroups.
\end{Thm}

Because of the use of the Fixed-point Lemma, one can only expect Bestvina-Handel theory to control the non-cyclic eigengroups of a
given automorphism of $F$. In fact, the cyclic ones will remain uncontrolled under the graph theoretic point of view.

\begin{Dfn}\label{defs5}
By Theorem~\ref{bh no connex}, every $\Phi\in Out(F)$ has finitely many non-cyclic conjugacy classes of eigengroups, say $k\geq 0$. A
\emph{set of representatives} for $\Phi$ is a set $\{ \phi_1, \ldots ,\phi_k \} \subset \Phi$ containing one and only one automorphism
in any isogredience class of $\Phi$ with non-cyclic fixed subgroup. Then, we denote by $\fix{\Phi}$ the corresponding set of conjugacy
classes of subgroups of $F$,
 $$
\fix{\Phi}=\{ [[\fix{\phi_1}]], \ldots ,[[\fix{\phi_k}]]\}
 $$
which, by Lemma~\ref{nc1}, does not depend on the particular set of representatives used. Note that, by construction, every non-cyclic
eigengroup of $\Phi$ is conjugate to one and only one of $\fix{\phi_1}, \ldots ,\fix{\phi_k}$. Also, by Theorem~\ref{bh no connex},
$r([[\fix{\phi_1}]])+\cdots +r([[\fix{\phi_k}]])\leq n+k-1$.
\end{Dfn}

\section{Free factor systems}\label{freefac}

Almost all the ideas in this section are essentially contained in subsection~2.6 of~\cite{BFH}. We restate and extend them here for
later use.

\begin{Dfn}\label{defs6}
Let $H,K\leq F$. We write $[[H]]\leq [[K]]$ (resp. $[[H]]\sqsubseteq [[K]]$) when $H^x\leq K^y$ (resp. $H^x$ is a free factor of
$K^y$) for some representatives $H^x\in [[H]]$ and $K^y\in [[K]]$. Note that, in this case, for every $x$ there exists $y$, and for
every $y$ there exists $x$, such that $H^x \leq K^y$ (resp. $H^x$ is a free factor of $K^y$). For example, the two unique conjugacy
classes of subgroups in $F$ with a single element (namely the {\it trivial} one $[[1]]=\{ 1\}$ and the {\it total} one $[[F]]=\{ F\}$)
are the extremal classes in these two partial orders, and $[[1]]\sqsubseteq [[F]]$.

A {\it subgroup system} of $F$ is a finite set of non-trivial conjugacy classes of subgroups of $F$. We will refer to the empty set as
the {\it trivial} subgroup system. Given two subgroup systems, $\mathcal{H}=\{ [[H_1]],\ldots ,[[H_r]]\}$, $\mathcal{K}=\{
[[K_1]],\ldots ,[[K_s]]\}$, we write $\mathcal{H}\leq \mathcal{K}$ if for every $i=1,\ldots ,r$ there exists $j=1,\ldots ,s$ such that
$[[H_i]]\leq [[K_j]]$.

A subgroup system $\mathcal{H}$ is called {\it finitely generated} when it only contains conjugacy classes of finitely generated
subgroups. In this case, the {\it complexity} of $\mathcal{H}$, denoted $cx(\mathcal{H})$, is defined as 0 if $\mathcal{H}$ is
trivial, and as the non-increasing sequence of positive integers obtained from the unordered list of ranks of the conjugacy classes in
$\mathcal{H}$, otherwise.

A subgroup system $\mathcal{H}=\{ [[H_1]],\ldots ,[[H_r]]\}$ is called a {\it free factor system} of $F$ when $\mathcal{H}=\emptyset$
or $H_1^{y_1} * \cdots *H_r^{y_r}$ is a free factor of $F$, for some choice $H_1^{y_1}\in [[H_1]], \ldots ,H_r^{y_r}\in [[H_r]]$. Note
that, in this event, this last condition is not satisfied in general for every conjugate. For example, $\{ [[\langle a\rangle]],
[[\langle b\rangle]]\}$ is a free factor system of $F=\langle a,b\rangle$ since $\langle a \rangle *\langle b\rangle$ is a free factor
of $F$, while $\langle a \rangle *\langle b^{ab}\rangle$ is not. Every free factor system except the total one, $\{[[F]]\}$, is called
\emph{proper}. That is, $\mathcal{H}$ is proper if, and only if, $H_i$ is a proper subgroup of $F$ for every $i=1,\ldots ,r$.

Every free factor system of $F$ is a finitely generated subgroup system and so it has a well defined complexity. It is clear that
there are only finitely many such complexities, namely the non-increasing sequences of positive integers adding up at most $n=r(F)$.
This finite set will be considered with the lexicographical order. For example, the trivial free factor system has the smallest
complexity, $cx(\emptyset )=0$, the total free factor system has the highest complexity, $cx(\{[[F]]\})=n$, and $cx(\{ [[\langle a_1
,a_2\rangle]], [[\langle a_3 \rangle]], \ldots ,[[\langle a_n \rangle]]\})=2,1,\stackrel{n-2}{\ldots},1$, where $n=r(F)\geq 3$ and $\{
a_1,\ldots ,a_n\}$ is a basis of $F$. We have $0<2,1,\stackrel{n-2}{\ldots},1<n$.

Given two subgroup systems of $F$, $$\mathcal{H}=\{ [[H_1]],\ldots ,[[H_r]]\} \, \,  \mbox{\rm and } \, \mathcal{K}=\{ [[K_1]],\ldots
,[[K_s]]\},$$ we write $\mathcal{H}\sqsubseteq \mathcal{K}$ if for every $i=1,\ldots ,r$ there exists $j=1,\ldots ,s$ such that
$[[H_i]]\sqsubseteq [[K_j]]$. Note that, if $\mathcal{H}$ and $\mathcal{K}$ are free factor systems then $j$ is uniquely determined by
$i$.

Let $H,K\leq F$ be finitely generated. By Proposition~2.1 in~\cite{N}, there are finitely many conjugacy classes of subgroups of $F$
of the type $H^x\cap K^y$, $x,y\in F$. This observation allows to define the {\it intersection} of two finitely generated conjugacy
classes $[[H]]$ and $[[K]]$ as
 $$
[[H]]\wedge [[K]]=\{ [[H^x \cap K^y]]\neq [[1]] \colon x,y\in F\}=\{ [[H\cap K^y]]\neq [[1]]\colon y\in F\}.
 $$
Note that $[[H]]\wedge [[K]]=\emptyset$ precisely when every conjugate of $H$ intersects trivially with every conjugate of $K$.
Analogously, if $\mathcal{H}=\{ [[H_1]],\ldots ,[[H_r]]\}$ and $\mathcal{K}=\{ [[K_1]],\ldots ,[[K_s]]\}$ are two finitely generated
subgroup systems, we define the {\it intersection} as
 $$
\mathcal{H}\wedge \mathcal{K}=\{ [[H_i \cap K_j^y]]\neq [[1]] : i=1,\ldots ,r,\,\,\, j=1,\ldots ,s,\,\,\, y\in F \}.
 $$
Clearly, $\mathcal{H}\wedge \mathcal{K}=\mathcal{K}\wedge \mathcal{H}\leq \mathcal{H}, \mathcal{K}$. In~\cite{BFH} Lemma~2.6.2, it is
proved that if $\mathcal{H}$ and $\mathcal{K}$ are free factor systems then so is $\mathcal{H}\wedge \mathcal{K}$.
\end{Dfn}

Automorphisms of $F$ act on the set of subgroups of $F$ in the natural way. Similarly, outer automorphisms act on the set of conjugacy
classes of subgroups of $F$. These observations point to the following definitions.

\begin{Dfn}\label{defs7}
Let $\Phi\in Out(F)$.

For every $H\leq F$, we write $[[H]]\Phi =[[H\phi]]$, where $\phi\in \Phi$. It is said that $[[H]]$ is {\it $\Phi$-invariant} when
$[[H]]\Phi =[[H]]$. In this case, for every $H^y\in [[H]]$ one can find $\phi \in \Phi$ such that $H^y$ is $\phi$-invariant.

In the same way, if $\mathcal{H}=\{ [[H_1]],\ldots ,[[H_r]]\}$ is a subgroup system of $F$, we define $\mathcal{H}\Phi$ as
$\mathcal{H}\Phi=\{ [[H_1]]\Phi,\ldots ,[[H_r]]\Phi\}$. And it is said that $\mathcal{H}$ is {\it $\Phi$-invariant} if $[[H_i]]$ is
$\Phi$-invariant for every $i=1,\ldots ,r$. Note that if $\mathcal{H}$ is a free factor system, then so is $\mathcal{H}\Phi$.

As noted in section~\ref{s BH}, the subgroup system $\{ [[\fix{\phi_1}]],$ $\ldots ,[[\fix{\phi_k}]]\}$ does not depend on the
particular set of representatives $\{ \phi_1, \ldots ,\phi_k\}$ chosen for $\Phi$. This is the {\it fixed subgroup system} of $\Phi$,
denoted $\fix{\Phi}$:
 $$
\begin{array}{rl}
\fix{\Phi} & =\{ [[\fix{\phi}]] : \phi\in \Phi, \,\,\, r(\fix{\phi})\geq 2 \} \\ & =\{ [[\fix{\phi_1}]], \ldots ,[[\fix{\phi_k}]]\}.
\end{array}
 $$
Finally, we observe that, for every $r\in \mathbb{Z}$ and $\phi\in \Phi$, $\fix{\phi}\leq \fix{\phi^r}$. Hence, $\fix{\Phi}\leq
\fix{\Phi^r}$.
\end{Dfn}

\begin{Emp}
An interesting example of a free factor system comes from the graph theoretic setting. Let $(Z,v_0,\eta)$ be a finite, connected,
marked graph, let $Z_0\leq Z$ be a subgraph, and denote by $Z_1, \ldots ,Z_s$ the non-contractible components of $Z_0$. For every
$i=1,\ldots ,s$, choose a vertex $v_i\in VZ_i$ and a path $q_i\in \pi Z(v_i,v_0)$, consider the inclusion $\pi Z_i (v_i,v_i)\to F$,
$x\mapsto (q_i^{-1}\cdot x\cdot q_i)\eta^{-1}$, and let $1\neq H_i\leq F$ be its image. It is clear that changing the chosen vertices
and paths, the $H_i$'s only change by conjugation. So, $\{ [[H_1]], \ldots ,[[H_s]]\}$ is a well-defined subgroup system of $F$,
determined by the subgraph $Z_0$, and denoted $\mathcal{H}(Z_0)$. By choosing maximal subtrees $T_i$ of $Z_i$, $i=1,\ldots, s$,
extending them to a maximal subtree $T$ of $Z$, and taking $q_i=T[v_i,v_0]$, we see that $H_1*\cdots *H_s$ is a free factor of $F$.
So, $\mathcal{H}(Z_0)$ is a free factor system of $F$. Furthermore, note that if $Z$ is a core graph and $Z_0$ is a proper subgraph of
$Z$, then $\mathcal{H}(Z_0)$ is proper.
\end{Emp}

Now, we prove the following two lemmas for later use. The last one is analogous to Lemma~\ref{iff}, but in the context of free factor
systems.

\begin{Lem}\label{inclusion ffs}
Let $F$ be a finitely generated free group and let $\mathcal{H}$ and $\mathcal{K}$ be two free factor systems of $F$ such that
$\mathcal{H}\leq \mathcal{K}$. Then, $\mathcal{H}\sqsubseteq \mathcal{K}$ and $cx(\mathcal{H})\leq cx(\mathcal{K})$, with strict
inequality except when $\mathcal{H}=\mathcal{K}$.
\end{Lem}

\demo Clearly, we may assume that $\mathcal{H}$ and $\mathcal{K}$ are both non-trivial. Write $\mathcal{H}=\{ [[H_1]],\ldots ,[[H_r]]
\}$ and $\mathcal{K}=\{ [[K_1]],\ldots ,[[K_s]] \}$, and define the map $\sigma \colon \{ 1,\ldots ,r\} \rightarrow \{ 1,\ldots ,s\}$
and choose elements $x_i\in F$ such that $H_i^{x_i}\leq K_{i\sigma}$. Since $H_i^{x_i}$ is a free factor of $F$, it also is a free
factor of $K_{i\sigma}$. Hence, $\mathcal{H}\sqsubseteq \mathcal{K}$. Furthermore, $\ast_{i\sigma=j} H_i^{x_i}$ is a free factor of
$K_j$, for every $j=1,\ldots ,s$.

As observed before, $i\sigma$ is uniquely determined by $i$. Hence, $cx(\mathcal{H})$ is obtained from $cx(\mathcal{K})$ by deleting
the entries corresponding to indices outside the image of $\sigma$, and replacing every other, say $t$, with a finite collection of
positive integers adding up at most $t$. Thus, $cx(\mathcal{H})\leq cx(\mathcal{K})$. And the inequality is strict except when
$\sigma$ is bijective and $H_i^{x_i}=K_{i\sigma}$ for all $i$, that is, except when $\mathcal{H}=\mathcal{K}$. \qed

\begin{Lem} \label{outiff}
Let $F$ be a finitely generated free group and let $\Phi\in Out(F)$. If there is a proper free factor system $\mathcal{K}$ of $F$ with
$\fix{\Phi}\leq \mathcal{K}$ then there exists a (proper) $\Phi$-invariant free factor system $\mathcal{H}$ of $F$ such that
$\fix{\Phi}\leq \mathcal{H}\sqsubseteq \mathcal{K}$.
\end{Lem}

\demo Let $\mathcal{H}=\{ [[H_1]],\ldots ,[[H_r]] \}$ be a (proper) free factor system satisfying $\fix{\Phi}\leq \mathcal{H}\leq
\mathcal{K}$ and having the smallest possible complexity. We only have to show that $\mathcal{H}$ is $\Phi$-invariant.

Let $\phi\in \Phi$.

For every $y\in F$ with $\fix{\phi \gamma_y}$ non-cyclic, there exist $i=1,\ldots ,r$ and an element $a\in F$ such that $\fix{\phi
\gamma_y}\leq H_i^a$. Also, $$\fix{\phi \gamma_y}\leq H_i^a \phi \gamma_y=(H_i \phi)^{(a\phi)y}$$ and so, $[[\fix{\phi \gamma_y}]]\leq
[[H_i \cap (H_i \phi )^b]]$ where $b=(a\phi)ya^{-1}\in F$. Hence, we have that $\fix{\Phi}\leq \mathcal{H}\wedge \mathcal{H}\Phi$.
Now, using Lemma~\ref{inclusion ffs} and the minimality of $cx(\mathcal{H})$, we have $\mathcal{H}\wedge \mathcal{H}\Phi =\mathcal{H}$
and so, $\mathcal{H}\leq \mathcal{H}\Phi$. But $cx(\mathcal{H})= cx(\mathcal{H}\Phi)$ so, again by Lemma~\ref{inclusion ffs},
$\mathcal{H}=\mathcal{H}\Phi$.

To complete the proof that $\mathcal{H}$ is $\Phi$-invariant, it remains to show that $\Phi$ does not permute the elements in
$\mathcal{H}$. Take $i=1,\ldots ,r$ and let $j$ and $c\in F$ be such that $H_i =(H_j \phi)^c$. By the minimality of $cx(\mathcal{H})$,
there exists $d\in F$ such that $1\neq \fix{\phi \gamma_d}\leq H_i =(H_j \phi)^c$. But also $\fix{\phi \gamma_d}\leq H_i \phi \gamma_d
=(H_i \phi )^d$. Thus, $1\neq \fix{\phi \gamma_d}\leq (H_j \phi)^c \cap (H_i \phi)^d$, which is only possible if $i=j$. This concludes
the proof. \qed

\section{Improved relative train tracks}\label{improved rtt}

Recently, Bestvina, Feighn and Handel published the paper~\cite{BFH} where, among other results, they made some improvements to
Bestvina-Handel theory. In the proof of our main result, we need two facts from this improved theory to deal with the exponential
case. This section is dedicated to them.

The next two results involve relative train tracks, which are equivalences of graphs with some special properties concerning a
filtration of the graph by invariant subgraphs (see~\cite{BH} for the precise definition). First, we restate the parts of
Theorem~5.1.5 and Lemma~5.1.7 in~\cite{BFH} that will be used later. In this reformulation, we only play attention to the top stratum,
however we restate the results, taking into account the technical differences used above: by subdividing at non-vertex fixed points
(see Proposition~III.2.3 in~\cite{DV}), we assume that the resulting relative train track has no such points; and adding the vertices
in the interior of the top stratum to $Z_0$, we assume that $VZ_0=VZ$ (obviously, these technical changes do not affect the other
properties).

\begin{Thm}[Bestvina-Feighn-Handel,~\cite{BFH}]\label{bfh}
Let $F$ be a non-tri\-vial fi\-ni\-te\-ly generated free group, and let $\Phi \in Out(F)$. Then, there exists a relative train track
$\beta$ of a non-empty, finite, connected, marked, core graph $(Z,v_0,\eta)$ such that $\Phi_{\beta}=\Phi^r$ for some $r\geq 1$, and
\begin{apartats}
\item all $\beta$-fixed points are vertices of $Z$, \item $Z$ has
a proper $\beta$-invariant subgraph $Z_0$ on which $\beta$ induces a self-map $\beta_0$, \item $\fix{\beta}$ has a basis $B$
containing (as a subgraph) a basis $B_0$ of $\fix{\beta_0}$, \item if $(\beta, Z,Z_0)$ is exponential then one of the following holds:
\begin{apartats}
\item[\upshape{(a)}] $B=B_0$, or \item[\upshape{(b)}] $B\setminus
B_0 =\{ p\}$ for some $p\in \pi Z$ with $\iota p=\tau p$ being an isolated vertex of $B_0$, or \item[\upshape{(c)}] $B\setminus B_0
=\{ p\}$ for some $p\in \pi Z$ with $\iota p\neq \tau p$, and one of $\iota p$ or $\tau p$ contained in a tree component of $Z_0$ (and
so, of $B_0$).
\end{apartats}
\end{apartats}
\end{Thm}

Note that case (b) corresponds to the geometric case of Theorem~5.1.5 in~\cite{BFH}, while (c) is the non-geometric one, for which
Lem\-ma~5.1.7 applies. However, we have restated those results by applying the cosmetic change whereby we have assumed $VZ_0=VZ$.

The following result, which is implicit in~\cite{BFH} but not explicitly stated, will be needed later. We thank M. Feighn for letting
us know about it and for helping us to extract a proof from~\cite{BFH}.

\begin{Thm}\label{exp}\label{cbfh}
Let $F$ be a non-trivial finitely generated free group, and let $\Phi \in Out(F)$. If every relative train track representing $\Phi$
has an exponential top stratum, then the same is true for all positive powers of $\Phi$.
\end{Thm}

\demo Assume that every relative train track representing $\Phi$ has an exponential top stratum, and that there exists $r\geq 2$ and a
relative train track $\beta\colon Z\to Z$, representing $\Phi^r$ and having a level top stratum (denote by $Z_0$ the maximal
$\beta$-invariant subgraph in the filtration). We will find a contradiction.

For every attracting lamination of $\Phi^r$, $\Omega \in \mathcal{L}(\Phi^r)$, (see Definition~3.1.5 in~\cite{BFH}), choose a
$\Omega$-generic line $\omega \in \Omega$, let $\hat{\omega}$ be the realization of $\omega$ in $Z$, and observe that, by
Lemma~3.1.10~(1) in~\cite{BFH}, $\hat{\omega}$ is a bi-infinite path running inside $Z_0$. So $\mathcal{L}(\Phi^r)$ is carried
(see~\cite{BFH}~page 532) by the proper $\Phi^r$-invariant free factor system $\mathcal{H}(Z_0)$. Let $\mathcal{F}$ be a proper free
factor system carrying $\mathcal{L}(\Phi^r)$ and with the smallest possible complexity. Since, by definition, $\mathcal{L}(\Phi^r)=
\mathcal{L}(\Phi)$, we see that $\mathcal{F}$ is $\Phi$-invariant.

Now apply Lemma~2.6.7 in~\cite{BFH} to obtain a relative train track map $\beta'\colon Z'\to Z'$ representing $\Phi$, with maximal
invariant subgraph $Z'_0$, and such that the free factor system $\mathcal{F} \sqsubset \mathcal{H}(Z'_0)$. By hypothesis, the top
stratum of $Z'$ is exponential. Utilising the argument in the proof of Lemma~3.1.13 in~\cite{BFH}, one can raise $\beta'$ to an
appropriate power $s\geq 1$, and obtain a relative train track representing $\Phi^s$ with maximal proper $\beta'^s$-invariant subgraph
$Z''_0$ containing $Z'_0$ and with the top stratum being exponential and aperiodic. In particular, $\mathcal{F} \sqsubset
\mathcal{H}(Z''_0)$. Now, consider the attracting lamination $\Omega\in \mathcal{L}(\Phi^s)$ associated to the top stratum of
$\beta'^s$ (see Definitions~3.1.12 in~\cite{BFH}). On the one hand, $\Omega$ is not carried by $\mathcal{H}(Z''_0)$ but, on the other
hand, $\Omega\in \mathcal{L}(\Phi^s)=\mathcal{L}(\Phi)$ is carried by $\mathcal{F}$. This contradiction concludes the proof. \qed

\section{The main argument}\label{s main theorem}

For all the section, let $F$ be a finitely generated free group of rank $n\geq 2$, and let $\Phi\in Out(F)$ such that
$r(\fix{\phi})\geq 2$ for some $\phi\in \Phi$, which is to say that $\fix{\Phi}$ is non-trivial in the sense of
Definitions~\ref{defs6} and~\ref{defs7}

Choose an equivalence $\beta\colon Z\to Z$ of a non-empty, finite, connected, marked, core graph $(Z,v_0,\eta)$, with all the
$\beta$-fixed points being vertices, and such that $\Phi_{\beta}=\Phi$ (see definitions~\ref{defs2} and~\ref{defs3}). Suppose also the
$\beta$ satisfies the following hypothesis (for example, the equivalences given in Theorems~\ref{threpr} and~\ref{bfh} do satisfy
them):

\begin{Hip}\label{hyp}
$Z$ has a proper $\beta$-invariant subgraph $Z_0$ on which $\beta$ induces a self-map $\beta_0$. Also, $\fix{\beta}$ has a basis $B$
having a subgraph $B_0$ which is a basis of $\fix{\beta_0}$, and such that either $B=B_0$ or $B\setminus B_0 =\{p\}$ for some
$\beta$-fixed path $p\in EB$.
\end{Hip}

Let us fix the following notation for the rest of the section.

\begin{Not} \label{notat}
(See section~\ref{s BH} for a review of groupoids, equivalences, representatives of outer automorphisms and relative train tracks).
Let $\{ \phi_1, \ldots, \phi_k\}$ be a set of representatives for $\Phi$; we have $k\geq 1$. Choose $p\in \pi Z(v_0\beta, v_0)$ such
that $\phi_{\beta, p,\eta}=\phi_1$, and consider the $\beta$-fixed vertex $v_1$ and the path $q_1\in \pi Z(v_1, v_0)$ with $q_1\beta
=q_1 \cdot p^{-1}$ given by the Fixed-point Lemma~\ref{fpl}. Changing the marking to $(Z,v_1,\eta')$, where $\eta'\colon F\to \pi
Z(v_1,v_1)$, $x\mapsto q_1\cdot x\eta \cdot q_1^{-1}$, we have $\phi_1 =\phi_{\beta,p,\eta}=\phi_{\beta,v_1,\eta'}=\eta'
\beta_{v_1}\eta'^{-1}$. From now on, $\beta$ will be considered as an equivalence of the marked graph $(Z,v_1,\eta')$ and $F$ will be
identified with $\pi Z(v_1,v_1)$ via $\eta'$ (so elements in $F$ are closed paths at $v_1$, with no further reference in the
notation). Then, $\phi_1 =\beta_{v_1}$ and $\fix{\phi_1}=\pi B(v_1,v_1)$. Furthermore, observe that if $v_1'$ is another $\beta$-fixed
vertex and $q\in \pi Z(v_1',v_1)$ is a $\beta$-fixed path then $(q\cdot q_1)\beta =(q\cdot q_1 )\cdot p^{-1}$; so $v_1$ can be
replaced by $v_1'$ and $\eta'$ by the new marking $\eta'' \colon F \to \pi_1 Z(v_1',v_1')$, $ x \mapsto q.q_1\cdot x\eta \cdot
(q.q_1)^{-1}$. Under this new marking we get as before that $\phi_1 =\beta_{v_1'}$ and $\fix{\phi_1}=\pi B(v_1',v_1')$. Hence, we may
interchange the roles of $v_1$ and $v_1'$.

By Proposition~\ref{corbij}, $B$ has exactly $k$ connected components with rank $\geq 2$. And, by the Fixed-point Lemma~\ref{fpl},
they can be labelled $B_1, \ldots ,B_k$ in such a way that $v_1\in VB_1$ and, for every $i=2,\ldots ,k$, there exists a vertex $v_i\in
VB_i$ and a path $q_i\in \pi Z(v_i,v_1)$ such that, $x\phi_i=q_i^{-1}\cdot (q_i \cdot x\cdot q_i^{-1})\beta_{v_i} \cdot q_i$ for every
$x\in F$, and $\fix{\phi_i}=q_i^{-1}\cdot \pi B_i (v_i,v_i)\cdot q_i$.
 $$
\xymatrix{
 \pi Z(v_1,v_1) \ar[dd]_{\phi_i} & \pi Z(v_i,v_i) \ar[l]_{\ \gamma_{q_i}} \ar[dd]^{\beta} \\
 \\ \pi Z(v_1,v_1) & \pi Z(v_i,v_i)
\ar[l]_{\ \gamma_{q_i}} }
 $$
Resetting $q_1=v_1$, the same equations are also valid for $i=1$. Additionally note that, changing $q_i$ into another path $q_i'\in
\pi Z(v_i,v_1)$ in the previous two equations, changes  $\phi_i$ to $\gamma_c^{-1}\phi_i \gamma_c \sim \phi_i$, where $c=q_i^{-1}\cdot
q_i'\in F$, and $\fix{\phi_i}$ gets right conjugated by $c$, $i=1,\ldots ,k$.

Renumbering if necessary, we can assume that either $B=B_0$, or $p$ is an edge of some component of $B$ with rank less than 2, or
$p\in EB_1$. Let $B_1'=B_1$ in the first two cases and $B_1'=B_1\setminus \{p\}$ in the last one. Let $B_i'=B_i$ for $i=2,\ldots ,k$.
In the case where $p\in EB_1$, the vertices $v_1$ and $\iota p$ both belong to $B_1$ so, we can apply the observation above to assume
that $v_1=\iota p$. We denote the terminal vertex of $p$ by $w$. We then distinguish two subcases: if $p$ does not separate $B_1$ then
$B_1'$ is connected and $r(B_1)=r(B_1')+1$. Otherwise, $B_1'$ has two components, say $B_{1,\iota}'$ for the one containing $v_1=\iota
p$, and $B_{1,\tau}'$ for the one containing $w=\tau p$; furthermore, $r(B_1)=r(B_{1,\iota}')+r(B_{1,\tau}')$.
\end{Not}

For a general $\beta$, this construction gives no information about the inclusions $\fix{\phi_i}\leq F$, since the edges of $B$ are
$\beta$-fixed elements in $\pi Z$ that can be arbitrarily complicated. However, hypothesis~\ref{hyp} will allow us to obtain some
algebraic information about those eigengroups.

\begin{Lem}\label{incpff}
Let us suppose that hypothesis~\ref{hyp} is satisfied. Using the notation above, if either
\begin{apartats}
\item $B=B_0$, or \item $B\setminus B_0=\{p\}$ and $p$ is an edge
of some component of $B$ with rank less than 2, or \item $B\setminus B_0=\{p\}$, $p\in EB_1$, $p$ separates $B_1$ and either
$B_{1,\iota}'$ or $B_{1,\tau}'$ is a tree,
\end{apartats}
then $\fix{\Phi}\leq \mathcal{H}(Z_0)$, which is a proper free factor system of $F$.
\end{Lem}

\demo Let $Z_1, \ldots ,Z_s$ be the components of $Z_0$ with rank at least 2. For every $i\in \{1,\ldots ,k\}$ let $j_i\in \{1,\ldots
,s\}$ be such that $v_i\in VZ_{j_i}$ (note that, in general, $j$ is neither one-to-one nor onto as a function of $i$, hence there is
no relation between $k$ and $s$). Then, $\pi B'_i\leq \pi Z_{j_i}$. In particular, for $i=1,\ldots ,k$, $$ q_i^{-1}\cdot \pi
B'_i(v_i,v_i)\cdot q_i \leq q_i^{-1}\cdot \pi Z_{j_i}(v_i,v_i)\cdot q_i. $$ Observe that, for $i\neq 1$, $B_i'=B_i$ and then the left
hand side of the above equation equals $\fix{\phi_i}$. Hence, $\{ [[\fix{\phi_2}]], \ldots ,[[\fix{\phi_k}]]\}\leq \mathcal{H}(Z_0)$.

It remains to show that also $\{[[\fix{\phi_1}]]\}\leq \mathcal{H}(Z_0)$. If (i) or (ii) occurs, then $B_1'=B_1$ and the argument
above also applies for $i=1$. Assume (iii). If $B_{1,\tau}'$ is a tree then $\pi B_1'(v_1,v_1)=\pi B_1 (v_1, v_1)$ and again the
argument above still works. So it only remains to consider the case where $B_{1,\iota}'$ is a tree (and $B_{1,\tau}'$ is not). We have
 $$
\fix{\phi_1}=\pi B_1(v_1,v_1)=p\cdot \pi B'_{1,\tau}(w,w)\cdot p^{-1}\leq p\cdot \pi Z_{j_1}(w,w)\cdot p^{-1}.
 $$
So, $[[\fix{\phi_1}]]$ is contained in the conjugacy class of subgroups of $F$ determined by the connected subgraph $Z_{j_1}$ of
$Z_0$. Thus $\fix{\Phi}\leq \mathcal{H}(Z_0)$.

Moreover, the free factor system $\mathcal{H}(Z_0)$ of $F$ is proper since $Z_0$ is a proper subgraph of $Z$ and $Z$ is a core graph.
\qed

This lemma will be applied to the equivalence given by Theorem~\ref{bfh}. Note that cases (iv)(a) and (iv)(b) there correspond to
cases (i) and (ii) here, respectively, and case (iv)(c) corresponds to (ii) or (iii) depending on whether or not the rank of the
component of $B$ containing $p$ is less than $2$.

To obtain the following result, it remains to closely analyse the path $p$, in the cases it exists. As will be seen, the case where
$(\beta, Z, Z_0)$ is level is the only one where $p$ plays a relevant role; additionally, by Theorem~\ref{threpr}(v), $p$ is a single
edge of $Z$ in this case.

\begin{Thm}\label{main}
Let $\Phi$ be an outer automorphism of a finitely generated free group $F$ such that $\fix{\Phi}$ is non-trivial; let $\phi\in \Phi$
with $r(\fix{\phi})\geq 2$. There exists a set of representatives $\{\phi_1,\ldots ,\phi_k\}$ for $\Phi$, $k\geq 1$, such that one of
the following holds:
\begin{apartats}
\item there is a proper $\Phi$-invariant free factor system
$\mathcal{F}$ such that $\fix{\Phi}\leq \mathcal{F}$,
\item there are two $\phi_1$-invariant subgroups $H,K\leq F$ such
that $F=H*K$; moreover, $\fix{\phi_1}=(H\cap \fix{\phi_1})* (K\cap \fix{\phi_1})$, $r(K\cap \fix{\phi_1})=1$ and, for each $i=2,\ldots
,k$, either $H$ is $\phi_i$-invariant and $\fix{\phi_i}\leq H$, or $K$ is $\phi_i$-invariant and $\fix{\phi_i}\leq K$,
\item there is a subgroup $H\leq F$ $\phi_i$-invariant for every
$i=1,\ldots ,k$, and non-trivial elements $y\in F$ and $h'\in H$, such that $F=H*\langle y\rangle$ and $y\phi_1 =h'y$; moreover,
$\fix{\phi_i}\leq H$ if $i\neq 1$, and $\fix{\phi_1}=(H\cap \fix{\phi_1})*\langle y^{-1}hy \rangle$ for some non-proper power $1\neq
h\in H$ with $h\phi_1 =h'hh'^{-1}$.
\end{apartats}
Furthermore, the choice can be made so that $\phi=\phi_i$ for some $i$.
\end{Thm}

\demo By Lemma~\ref{imagey}, it is sufficient to prove the Theorem with the following weaker statement in place of (iii):
\begin{apartats}
\item[\upshape{(iii')}] \emph{there is a subgroup $H\leq F$
$\phi_i$-invariant for every $i=1,\ldots ,k$, and an element $1\neq y\in F$, such that $F=H*\langle y\rangle$; moreover,
$\fix{\phi_i}\leq H$ if $i\neq 1$, and $\fix{\phi_1}=(H\cap \fix{\phi_1})*\langle y^{-1}hy \rangle$ for some $1\neq h\in H$. }
\end{apartats}
The proof then works by induction on $n$, the rank of $F$. The result is vacuous for $n=0,1$. For $n=2$, the result is true using
Theorem~\ref{maxrank} (in this case, $k=1$, (i) cannot happen, and the two possibilities in Collins-Turner Theorem are special cases
of (ii) and (iii') here). So, we can assume $n\geq 3$, and that the result is true for free groups of smaller rank.

Note that the result involves both a choice of representatives for $\Phi$, and a choice of basis for $F$. Given such a basis
$\{x_1,\ldots ,x_n\}$ and a set of representatives $\{ \phi_1, \ldots ,\phi_k \}$, take an element $c\in F$ and consider the new basis
$\{x_1^c,\ldots ,x_n^c\}$ and the new set of representatives $\{ \phi'_1, \ldots ,\phi'_k\}$, where $\phi'_i =\gamma_{c}^{-1}\phi_i
\gamma_c\sim \phi_i$. This new set of automorphisms acts on the new basis just as the old one acted on the original basis. Hence, up
to a change of basis, we may choose one of the representatives to be our favorite $\phi \in \Phi$, as long as $r(\fix{\phi})\geq 2$.
Thus the last statement of the theorem follows as long as we can find some set of representatives for $\Phi$ satisfying the theorem.

Let $(\beta, Z, Z_0)$ be a minimal representative of $\Phi$ (the existence is ensured by Corollary~IV.1.2 of~\cite{DV}). Let us note
that, by Theorem~\ref{threpr}, $\beta$ satisfies hypothesis~\ref{hyp}. We shall henceforth use the Notation~\ref{notat}
 applied to $\beta$. We recall that we chose $\{\phi_1,\ldots ,\phi_k\}$ to be an arbitrary set of representatives for $\Phi$. Hence,
in order to prove the result, we are still free to change it at our convenience.

By Theorem~\ref{threpr}~(iv), $(\beta,Z,Z_0)$ is not null; hence, it is either level or exponential.

Suppose now that $(\beta,Z,Z_0)$ is exponential. By the minimality of $PF(\beta/Z_0)$, all the representatives of $\Phi$ are also
exponential. Consider the integer $r\geq 1$ and the representative $(\beta',Z',Z'_0)$ of $\Phi^r$ given by Theorem~\ref{bfh}.
Theorem~\ref{cbfh} implies that $(\beta',Z',Z'_0)$ is exponential and hence, \ref{bfh}~(iv) applies and we have (a), (b) or (c). Thus,
$\beta'$ is in the situation (i), (ii) or (iii) of Lemma~\ref{incpff} (see the subsequent comment there). Furthermore, $\beta'$
satisfies hypothesis~\ref{hyp}. Thus, we conclude that $\fix{\Phi^r}\leq \mathcal{H}(Z_0')$, a proper free factor system of $F$.
Finally, using the fact $\fix{\Phi}\leq \fix{\Phi^r}$ and Lemma~\ref{outiff}, we end up in case (i) of the theorem.

So, we can assume that $(\beta,Z,Z_0)$ is level. By Theorem~\ref{threpr}~(v), either $B=B_0$ or $B\setminus B_0 =\{ e \}$, where $e$
is a $\beta$-fixed edge such that $EZ\setminus EZ_0=\{ e\}$. Using Lemmas~\ref{incpff} and~\ref{outiff}, the first possibility and
some particular cases of the second one immediately leads us to case (i) of the theorem. So, we can assume that $B\setminus
B_0=\{e\}=EZ\setminus EZ_0$, $e\in EB_1$ and either $e$ does not separate $B_1$, or it does but into two non-tree components.

Now we distinguish two cases, depending on whether $e$ does or does not separate $Z$ (caution, not $B$ !).

\medskip

\noindent \emph{Case 1}: $e$ separates $Z$ (into two components).

Recall that we are using the notation established in~\ref{notat}. Let $X$ be the component of $Z_0=Z\setminus \{e\}$ containing
$v_1=\iota e$, and $Y$ be the one containing $w=\tau e$. Since $Z$ is a core graph and $e\beta =e$, both $X$ and $Y$ are
$\beta$-invariant subgraphs of $Z$ and neither is a tree. It is clear that $B_1$ is the unique component of $B$ with rank $\geq 2$
that has vertices both in $X$ and $Y$. So, renumbering if necessary, we can write $\{ B_1, \ldots ,B_k \}=\{B_1\} \cup \{B_2, \ldots
,B_l\}\cup \{ B_{l+1}, \ldots ,B_k\}$, where $\pi B_i\leq \pi X$ for $i=2,\ldots ,l$ and $\pi B_i\leq \pi Y$ for $i=l+1,\ldots ,k$,
$1\leq l\leq k$. Now, changing $\phi_2,\ldots ,\phi_k$ to appropriate isogredient automorphisms if necessary, we can assume that, for
$i=2,\ldots ,l$, $q_i\in \pi X(v_i,v_1)$ and for $i=l+1,\ldots ,k$, $q_i \in \pi Y(v_i, w)\cdot e^{-1}$. Moreover, the situation
forces $B_1'$ to be disconnected and recall that, in this case, we were assuming that neither $B'_{1,\iota}$ nor $B'_{1,\tau}$ is a
tree.

Let $L=\pi X(v_1,v_1)\leq F$ and $M=e\cdot \pi Y(w,w)\cdot e^{-1}\leq F$. It is straightforward to verify that $L,M\neq 1$ are both
$\phi_1$-invariant and that $F=L*M$. Furthermore,
 $$
\begin{array}{rl}
\fix{\phi_1} & =\pi B_1(v_1,v_1) \\ & =\pi B'_{1,\iota}(v_1,v_1)*(e\cdot \pi B'_{1,\tau}(w,w)\cdot e^{-1}) \\ & =(L\cap
\fix{\phi_1})*(M\cap \fix{\phi_1}),
\end{array}
 $$
where $L\cap \fix{\phi_1}\neq 1$ and $M\cap \fix{\phi_1}\neq 1$. From the equalities $$x\phi_i=q_i^{-1}\cdot (q_i \cdot x\cdot
q_i^{-1})\beta \cdot q_i \ \ \mbox{\rm and } \ \fix{\phi_i}=q_i^{-1}\cdot \pi B_i (v_i,v_i)\cdot q_i,$$ and from the
$\beta$-invariance of $X$ and $Y$, it is easy to see that, for $i=2,\ldots ,l$, $L$ is $\phi_i$-invariant and $\fix{\phi_i}\leq L$
and, for $i=l+1,\ldots ,k$, $M$ is $\phi_i$-invariant and $\fix{\phi_i}\leq M$.

If $r(L\cap \fix{\phi_1})=1$ then the set of representatives $\{ \phi_1 ,\ldots ,\phi_k\}$ for $\Phi$ satisfies case (ii) of the
theorem with $H=M$ and $K=L$. So, we can assume that $r(L\cap \fix{\phi_1})\geq 2$. In particular, $r(L)\geq 2$.

Let us consider now the restriction of $\beta$ to $X$, $\beta_{_X}\colon X\to X$. Since $L$ is a free factor of $F$, the restriction
of $\phi_1$ to $L$ is an automorphism of $L$ and hence, $\beta_{_X}$ is an equivalence of $X$. Let $\Phi_{_L}=\Phi_{\beta_{_X}}\in
Out(L)$. Alternatively, $\Phi_{_L}$ can also be defined by simultaneously restricting $\phi_1, \ldots ,\phi_l$ to $L$, and using
Lemma~\ref{malnormal} to see that all these restrictions belong to the same coset of $Inn(L)$. Although the situation with
$\beta_{_X}$ is not the same as it was for $\beta$ (for example, $X$ is not in general a core graph, and we know nothing in general
about maximal $\beta_{_X}$-invariant subgraphs of $X$), we can still apply Proposition~\ref{corbij} to $\beta_{_X}$ and $B'_{1,\iota}
\cup B_2 \cup \cdots \cup B_l$ (a basis for $\fix{\beta_{_X}}$), and deduce that $\{\phi_{1_L}, \ldots ,\phi_{l_L}\}$ is a set of
representatives for $\Phi_{_L}$, and that
 $$
\fix{\Phi_{_L}}= \{[[L\cap \fix{\phi_1}]], [[\fix{\phi_2}]],\ldots ,[[\fix{\phi_l}]]\}.
 $$
(recall that $\fix{\phi_2}, \ldots ,\fix{\phi_l}$ are all contained in $L$). Now, it is time to apply the inductive hypothesis to
$\Phi_{_L}$. Thus, we may find a set of representatives $\{\varphi_1, \ldots ,\varphi_l\}$ for $\Phi_{_L}$ such that
$\varphi_1=\phi_{1_L}, \varphi_2 \sim \phi_{2_L}, \ldots ,\varphi_l \sim \phi_{l_L}$, and their fixed subgroups satisfy the conclusion
of the theorem. Also, changing the automorphisms $\phi_2, \ldots ,\phi_l$ to isogredient ones (i.e. changing the paths $q_i$
appropriately in $\pi X(v_i,v_1)$, $i=2,\ldots ,l$) if necessary, we can assume that $\varphi_2 =\phi_{2_L}, \ldots ,\varphi_l
=\phi_{l_L}$. So, the situation now is the same as before the application of the inductive hypothesis, but with the extra information
that the set of representatives $\{\phi_{1_L}, \ldots ,\phi_{l_L}\}$ for $\Phi_{_L}$ satisfy (i), or (ii), or (iii') of the theorem.
Let us distinguish these three subcases.

\medskip

\emph{Subcase 1.1}: $\fix{\Phi_{_L}}\leq \mathcal{L}$ for some proper $\Phi_{_L}$-invariant free factor system $\mathcal{L}=\{
[[L_1]], \ldots ,[[L_t]]\}$ of $L$. Choose notation such that $L\cap \fix{\phi_1}\leq L_1$, $\fix{\phi_2}\leq L_{j_2}^{x_2},\ldots
,\fix{\phi_l}\leq L_{j_l}^{x_l}$ and $L_1*L_2^{y_2}*\cdots *L_t^{y_t}$ is a free factor of $L$, where $x_2,\ldots ,x_l, y_2, \ldots
,y_t \in L$. Then, $(L_1*M)*L_2^{y_2}*\cdots *L_t^{y_t}$ is a free factor of $F=L*M$, which means that $\mathcal{F}=\{[[L_1 *M]],
[[L_2]], \ldots ,[[L_t]]\}$ is a proper free factor system of $F$. And it is clear that $\fix{\Phi}\leq \mathcal{F}$. So, using
Lemma~\ref{outiff}, we end up in case (i) of the theorem.

\medskip

\emph{Subcase 1.2}: writing $\{\phi_{1_L}, \ldots ,\phi_{l_L}\}=\{ \psi_1, \ldots ,\psi_{l}\}$, there are two $\psi_1$-invariant
subgroups $C,D\leq L$ such that $L=C*D$; moreover, $$\fix{\psi_1}=(C\cap \fix{\psi_1})*(D\cap \fix{\psi_1}),$$ where $r(D\cap
\fix{\psi_1})=1$ and for each $j=2,\ldots ,l$, either $C$ is $\psi_j$-invariant and $\fix{\psi_j}\leq C$, or $D$ is $\psi_j$-invariant
and $\fix{\psi_j}\leq D$. Let us consider now the different possibilities, depending on which of the $\phi_{1_L},\ldots ,\phi_{l_L}$
is equal to $\psi_1$, and on the relationship between $\phi_{1_L}$ and the subgroups $C$ and $D$.

First suppose that $\psi_1 =\phi_{1_L}$. Let $H=C*M$ and $K=D$. We have $$F=L*M=C*D*M=H*K.$$ Since $C$, $D$ and $M$ are
$\phi_1$-invariant, so are $H$ and $K$. Also,
 $$
\begin{array}{rl}
\fix{\phi_1} & =(L\cap \fix{\phi_1})*(M\cap \fix{\phi_1}) \\ & =(C\cap \fix{\phi_1})*(D\cap \fix{\phi_1})*(M\cap \fix{\phi_1}) \\ &
=(H\cap \fix{\phi_1})*(K\cap \fix{\phi_1}),
\end{array}
 $$
where $r(K\cap \fix{\phi_1})=1$. For every $i=2,\ldots ,l$, we know that either $C$ is $\phi_i$-invariant and $\fix{\phi_i}\leq C$, or
$D$ is $\phi_i$-invariant and $\fix{\phi_i}\leq D$. In the first case $H$ is also $\phi_i$-invariant by an application of
Lemma~\ref{malnormal} to $\phi_1$ and $\phi_i$, and $\fix{\phi_i}\leq C\leq H$. And in the second case, we directly have that $K=D$ is
$\phi_i$-invariant and $\fix{\phi_i}\leq K$. On the other hand, for every $i=l+1,\ldots ,k$, we know that $M$ is $\phi_i$-invariant
(and then $H$, again by Lemma~\ref{malnormal}) and additionally  $\fix{\phi_i}\leq M\leq H$. Hence, the set of representatives $\{
\phi_1, \ldots ,\phi_k\}$ for $\Phi$ satisfies case (ii) of the theorem.

Now suppose that $\psi_1 \neq \phi_{1_L}$ and that $C$ is $\phi_1$-invariant and $L\cap \fix{\phi_1}\leq C$. Renumbering if necessary,
we can assume that $\psi_1 =\phi_{2_L}$. As before, we let $H=C*M$ and $K=D$ (and hence $F=H*K$). Since $C$ and $M$ are
$\phi_1$-invariant, so is $H$. Also, $\fix{\phi_1}=(L\cap \fix{\phi_1})*(M\cap \fix{\phi_1})\leq H$. For $\phi_2$ we have
 $$
\begin{array}{rl}
\fix{\phi_2} & =L\cap \fix{\phi_2} \\ & =(C\cap \fix{\phi_2})*(D\cap \fix{\phi_2}) \\ & =(H\cap \fix{\phi_2})*(K\cap \fix{\phi_2}),
\end{array}
 $$
where $r(K\cap \fix{\phi_2})=1$. Consequently, $H\cap \fix{\phi_2}\neq 1$ and, since $H$ is $\phi_1$-invariant, Lemma~\ref{malnormal}
implies that $H$ is also $\phi_2$-invariant; and we know that $K=D$ is also $\phi_2$-invariant. Furthermore, and using the same
arguments as in the previous paragraph, for every $i=3,\ldots ,k$, either $H$ is $\phi_i$-invariant and $\fix{\phi_i}\leq H$, or $K$
is $\phi_i$-invariant and $\fix{\phi_i}\leq K$. Hence, the (ordered) set of representatives $\{ \phi_2, \phi_1, \phi_3, \ldots
,\phi_k\}$ for $\Phi$ satisfies case (ii) of the theorem, with $\phi_2$ playing the distinguished role.

Finally, suppose that $\psi_1 \neq \phi_{1_L}$ and that $D$ is $\phi_1$-invariant and $L\cap \fix{\phi_1}\leq D$. Exactly the same
argument in the previous paragraph, interchanging $H$ with $K$ and $C$ with $D$, shows that $\{ \phi_2, \phi_1, \phi_3, \ldots
,\phi_k\}$ is a set of representatives for $\Phi$ satisfying case (ii) of the theorem with $H=C$ and $K=D*M$.

\medskip

\emph{Subcase 1.3}: there is a subgroup $K\leq L$ $\phi_i$-invariant for every $i=1,\ldots ,l$, and an element $1\neq y\in L$, such
that $L=K*\langle y\rangle$ (so, $r(K)=r(L)-1\geq 1$); moreover, there is an index $j\in \{1,\ldots ,l\}$ such that $L\cap
\fix{\phi_i}\leq K$ for $i\neq j$, and $L\cap \fix{\phi_j}= (K\cap \fix{\phi_j})*\langle y^{-1}hy\rangle$ for some $1\neq h\in K$.
Renumbering if necessary, we can assume that either $j=1$ or $j=2$.

Let $H=K*M$. We have $F=L*M=K*\langle y\rangle *M=H*\langle y\rangle$. We already know that $K$ and $M$ are $\phi_1$-invariant and
thus so is $H$. Also, for every $i=2,\ldots ,k$, $H\cap \fix{\phi_i}\neq 1$ and hence, by Lemma~\ref{malnormal}, $H$ is
$\phi_i$-invariant.

Suppose $j=1$. For every $i=2,\ldots ,l$, $\fix{\phi_i}=L\cap \fix{\phi_i}\leq K\leq H$, and for $i=l+1,\ldots ,k$, $\fix{\phi_i}\leq
M\leq H$. Also,
 $$
\begin{array}{rl}
\fix{\phi_1} & =(L\cap \fix{\phi_1})*(M\cap \fix{\phi_1}) \\ & =(K\cap \fix{\phi_1})*\langle y^{-1}hy\rangle *(M\cap \fix{\phi_1}) \\
& =(H\cap \fix{\phi_1})*\langle y^{-1}hy\rangle.
\end{array}
 $$
Thus, the set of representatives $\{\phi_1, \ldots ,\phi_k\}$ for $\Phi$ satisfies case (iii') of the theorem.

Suppose $j=2$. We have $\fix{\phi_1}=(L\cap \fix{\phi_1})*(M\cap \fix{\phi_1})\leq K*M=H$, and also $\fix{\phi_i}\leq H$ for every
$i=3,\ldots ,k$, just as in the previous paragraph. Furthermore,
 $$
\fix{\phi_2}=L\cap \fix{\phi_2}=(K\cap \fix{\phi_2})*\langle y^{-1}hy\rangle =(H\cap \fix{\phi_2})*\langle y^{-1}hy\rangle.
 $$
Thus, the set of representatives $\{\phi_2, \phi_1, \phi_3, \ldots ,\phi_k\}$ for $\Phi$ satisfies case (iii') of the theorem.

This concludes case 1 of the proof.

\medskip

\noindent \emph{Case 2}: $e$ does not separate $Z$.

In this case, $Z_0$ is connected. Consider $L=\pi Z_0(v_1,v_1)$, a free factor of $F$ with rank $n-1$. Changing $\phi_2,\ldots
,\phi_k$ to appropriate isogredient automorphisms if necessary, we can assume that the paths $q_i$ do not cross the edge $e$, i.e.
$q_i\in \pi Z_0(v_i,v_1)$. So, from the expression $x\phi_i=q_i^{-1}\cdot (q_i \cdot x\cdot q_i^{-1})\beta \cdot q_i$ and the
$\beta$-invariance of $Z_0$, we see that $L$ is $\phi_i$-invariant for every $i=1,\ldots ,k$. Moreover, for $i\neq 1$, the equality
 $$
\fix{\phi_i}=q_i^{-1}\cdot \pi B_i (v_i,v_i)\cdot q_i =q_i^{-1}\cdot \pi B_i' (v_i,v_i)\cdot q_i
 $$
tells us that $\fix{\phi_i}\leq L$.

Suppose that $e$ does not separate $B_1$ (and so $B_1'$ is connected). Choose an $r\in \pi B_1'(w,v_1)$ and let $z=e\cdot r\in
\fix{\phi_1}$, a non-trivial element of $F$. Since $z$ is a path crossing $e$ only once, it is clear that $F=L*\langle z\rangle$.
Furthermore,
 $$
\fix{\phi_1}=\pi B_1(v_1,v_1)=\pi B_1'(v_1,v_1)*\langle e\cdot r\rangle =(L\cap \fix{\phi_1})*\langle z\rangle.
 $$
Thus, taking $H=L$ and $K=\langle z\rangle$, we are in case (ii) of the theorem.

So, we can assume that $e$ separates $B_1$ into $B_{1,\iota}'$ and $B_{1,\tau}'$. As noted above, by using Lemmas~\ref{incpff}
and~\ref{outiff}, we have reduced to the case where neither $B_{1,\iota}'$ nor $B_{1,\tau}'$ is a tree. Let $r\in \pi Z_0(w, v_1)$.
The following arguments will work for every such path, up to a certain point in the proof, when we will choose a specific one to work
with. Let $1\neq z=e\cdot r\in F$. As before, $F=L*\langle z\rangle$ (the difference now is that $z$ is not fixed by $\phi_1$).
Consider the automorphism $\phi'_1\in \Phi$ given by the $\beta$-fixed vertex $w$ and the path $r\in \pi Z(w,v_1)$, i.e. $x\phi'_1
=r^{-1}\cdot (r \cdot x\cdot r^{-1})\beta \cdot r$, $x\in F$ (and note that $\fix{\phi'_i}=r^{-1}\cdot \pi B_1(w,w)\cdot r$). Observe
that $\fix{\phi'_1}=(\fix{\phi_1})^z$, that $\phi_1 \sim \phi'_1 =\gamma_{z}^{-1}\phi_1 \gamma_z$. Also note that as $r$ ranges over
all paths in $\pi Z_0(w,v_1)$, $\phi'_1$ ranges over all the automorphisms of the form $\gamma_{y}^{-1}\phi_1 \gamma_y$, $y\in L$.
Since $r\in \pi Z_0$, $L$ is also $\phi'_1$-invariant. We have
 $$
\begin{array}{rl}
\fix{\phi_1}& =\pi B_1(v_1,v_1) \\ & =\pi B_{1,\iota}'(v_1,v_1) *e\cdot \pi B_{1,\tau}'(w,w)\cdot e^{-1} \\ & =(L\cap \fix{\phi_1})*
(r^{-1}\cdot \pi B_{1,\tau}'(w,w)\cdot r)^{r^{-1}\cdot e^{-1}} \\ & =(L\cap \fix{\phi_1}) *(L\cap \fix{\phi'_1})^{z^{-1}}.
\end{array}
 $$
If $r(\pi B_{1,\tau}'(w,w))=1$ then $L\cap \fix{\phi'_1}=\langle h\rangle$ for some $h\in L$. Hence, $$\fix{\phi_1}=(L\cap
\fix{\phi_1})*\langle zhz^{-1} \rangle$$ allowing us to conclude that we are in case (iii') of the theorem with $H=L$ and $y=z^{-1}$.
Thus we can assume that $r(\pi B_{1,\tau}'(w,w))=r(L\cap \fix{\phi'_1})\geq 2$; in particular, $n-1=r(L)\geq 2$. Since
 $$
\fix{\phi'_1}=(\fix{\phi_1})^z=(L\cap \fix{\phi_1})^z*(L\cap \fix{\phi'_1}),
 $$
the same argument works interchanging $\phi_1$ with $\phi'_1$ and $z$ with $z^{-1}$. So, without loss of generality, we may assume
$r(\pi B_{1,\iota}'(v_1,v_1))=r(L\cap \fix{\phi_1})\geq 2$. With these assumptions, the components of $B_0$ with rank $\geq 2$ are
precisely $B_{1,\iota}', B_{1,\tau}', B_2, \ldots ,B_k$.

Let us consider now the restriction of $\beta$ to $Z_0$, $\beta_0\colon Z_0\to Z_0$. As we argued in case 1 of the proof, $\beta_0$ is
an equivalence of $Z_0$, which induces an outer automorphism of $L$ which we denote by
 $\Phi_{_L}=\Phi_{\beta_0}$. Then, by
Proposition~\ref{corbij} applied to $\beta_0$ and $B_0$, we deduce that $\{\phi_{1_L}, \phi'_{1_L}, \phi_{2_L}, \ldots ,\phi_{k_L}\}$
is a set of representatives for $\Phi_{_L}$ and
 $$
\fix{\Phi_{_L}}= \{[[L\cap \fix{\phi_1}]], [[L\cap \fix{\phi'_1}]],[[\fix{\phi_{2}}]],\ldots ,[[\fix{\phi_{k}}]]\}
 $$
(recall that $\fix{\phi_2}, \ldots ,\fix{\phi_k}$ are all contained in $L$). By the inductive hypothesis applied to $\Phi_{_L}$, we
may find a set of representatives, $\{\varphi_1, \varphi'_1, \varphi_2, \ldots ,\varphi_k \}$, for $\Phi_{_L}$, such that
$\varphi_1=\phi_{1_L}, \varphi'_1 \sim \phi'_{1_L}, \varphi_2 \sim \phi_{2_L}, \ldots ,\varphi_k \sim \phi_{k_L}$, and their fixed
subgroups satisfy the conclusion of the theorem. Recall that we are still free to change the choice of $r\in \pi Z_0(w,v_1)$. Doing
this appropriately, we may assume that $\varphi_1'=\phi'_{1_L}$. Also, changing the automorphisms $\phi_2, \ldots, \phi_k$ to
isogredient ones, if necessary, we may also assume that $\varphi_2 =\phi_{2_L}, \ldots ,\varphi_k =\phi_{k_L}$. So, the situation now
is the same as before the application of the inductive hypothesis, but with the extra information that the set of representatives
$\{\phi_{1_L}, \phi'_{1_L}, \phi_{2_L},\ldots ,\phi_{k_L}\}$ for $\Phi_{_L}\in Out(L)$ satisfy (i), or (ii), or (iii') of the theorem.
Let us distinguish these three subcases.

\medskip

\noindent \emph{Subcase 2.1}: $\fix{\Phi_{_L}}\leq \mathcal{L}$ for some proper $\Phi_{_L}$-invariant free factor system
$\mathcal{L}=\{ [[L_1]], \ldots ,[[L_t]]\}$ of $L$. Choose notation such that $L\cap\fix{\phi_1}\leq L_1$, $L\cap \fix{\phi'_1}\leq
L_{j_1}^{x_1}$, $\fix{\phi_{2}}\leq L_{j_2}^{x_2},\ldots ,\fix{\phi_{k}}\leq L_{j_k}^{x_k}$ and $L_1*L_2^{y_2}*\cdots *L_t^{y_t}$ is a
free factor of $L$, where $x_2,\ldots ,x_k, y_2, \ldots ,y_t\in L$. Then, $L_1*L_2^{y_2}*\cdots *L_t^{y_t}*\langle z\rangle$ is a free
factor of $F=L*\langle z\rangle$.

Suppose that $j_1=1$, and consider the proper free factor system $$\mathcal{F}=\{[[L_1*\langle z\rangle]], [[L_2]], \ldots
,[[L_t]]\}.$$ By changing the choice of $r\in \pi Z_0 (w,v_1)$ (and hence $z$ and $\phi'_1$ which depend on $r$) if necessary, we can
assume $x_1=1$. Now, $$\fix{\phi_1}=(L\cap \fix{\phi_1})*(L\cap \fix{\phi'_1})^{z^{-1}}\leq L_1*\langle z\rangle$$ and thus
$\fix{\Phi}\leq \mathcal{F}$. So, by Lemma~\ref{outiff}, we end up in case (i) of the theorem.

Otherwise, suppose that $j_1\neq 1$, say $j_1=2$, and note that $x_1$ can be taken to be trivial. Since $y_2\in L$, we have that
$L_1*L_2^{y_2}*\cdots *L_t^{y_t}*\langle zy_2\rangle$ and then $L_1*L_2^{z_{\,}^{-1}}*L_3^{y_3}*\cdots *L_t^{y_t}*\langle zy_2\rangle$
are also free factors of $F$. Hence, $$\mathcal{F}=\{[[L_1*L_2^{z_{\,}^{-1}}]], [[L_3]], \ldots ,[[L_t]]\}$$ is a proper free factor
system of $F$, and $\fix{\Phi}\leq \mathcal{F}$. Again, by Lemma~\ref{outiff}, we are in case (i) of the theorem.

\medskip

\noindent \emph{Subcase 2.2}: writing $\{\phi_{1_L}, \phi'_{1_L}, \phi_{2_L},\ldots ,\phi_{k_L}\}=\{ \psi_1, \ldots ,\psi_{k+1}\}$,
there are two $\psi_1$-invariant subgroups $C,D\leq L$ such that $L=C*D$; moreover, $$\fix{\psi_1}=(C\cap \fix{\psi_1})*(D\cap
\fix{\psi_1}),$$ where $r(D\cap \fix{\psi_1})=1$ and for each $j=2,\ldots ,k+1$, either $C$ is $\psi_j$-invariant and
$\fix{\psi_j}\leq C$, or $D$ is $\psi_j$-invariant and $\fix{\psi_j}\leq D$.

Recall that $\phi'_1=\gamma_{z}^{-1}\phi_1 \gamma_z =\phi_1 \gamma_g$ where $g=(z\phi_1)^{-1}z\in F$. Recall also that
$$\fix{\phi_1}=(L\cap \fix{\phi_1})*(L\cap \fix{\phi'_1})^{z^{-1}} \ \  \mbox{\rm and } \ \fix{\phi'_1}=(\fix{\phi_1})^z.$$ Let us
consider now the different possibilities, depending on which of the automorphisms $\phi_{1_L},\phi'_{1_L},\phi_{2_L},\ldots
,\phi_{k_L}$ is equal to $\psi_1$, and on the relationship between $\phi_{1_L} ,\phi'_{1_L}$ and $C,D$.

Suppose that $\psi_1 \neq \phi_{1_L}, \phi'_{1_L}$ and that $C$ is both $\phi_1$- and $\phi'_1$-invariant and $L\cap \fix{\phi_1},
L\cap \fix{\phi'_1}\leq C$. Renumbering if necessary, we can assume that $\psi_1 =\phi_{2_L}$. Let $H=C*\langle z\rangle$ and $K=D$.
We have $$F=L*\langle z\rangle =C*D*\langle z\rangle =H*K.$$ Since $C$ is $\phi_1$-invariant and $C\cap \fix{\phi'_1}\neq 1$,
Lemma~\ref{malnormal} implies that $g\in C$. So, $z\phi_1\in C*\langle z\rangle$ and $H$ is $\phi_1$-invariant. Also,
$$\fix{\phi_1}=(L\cap \fix{\phi_1})*(L\cap \fix{\phi'_1})^{z^{-1}}\leq C*\langle z\rangle =H.$$ For $\phi_2$ we know that
$C\cap\fix{\phi_2}=H\cap L\cap \fix{\phi_2}=H\cap \fix{\phi_2}$ and $$\fix{\phi_2}=(C\cap \fix{\phi_2})*(D\cap \fix{\phi_2})= (H\cap
\fix{\phi_2})*(K\cap \fix{\phi_2}),$$ where $r(K\cap \fix{\phi_2})=1$ (and so $H\cap \fix{\phi_2}\neq 1$). Additionally, $H$ is
$\phi_1$-invariant and intersects $\fix{\phi_2}$ non-trivially. Thus, Lemma~\ref{malnormal} implies that $H$ is $\phi_2$-invariant,
while $K=D$ is $\phi_2$-invariant by hypothesis. Finally, for each $i=3,\ldots ,k$, we know that either $C$ is $\phi_i$-invariant and
$\fix{\phi_i}=L\cap \fix{\phi_i}\leq C$, or $D$ is $\phi_i$-invariant and $\fix{\phi_i}\leq D$. In the first case $H$ is also
$\phi_i$-invariant by another application of Lemma~\ref{malnormal} to $\phi_1$ and $\phi_i$, and $\fix{\phi_i}\leq C\leq H$. And in
the second, by hypothesis, $K=D$ is $\phi_i$-invariant and $\fix{\phi_i}\leq K$. Hence, the set of representatives $\{ \phi_2, \phi_1,
\phi_3, \ldots ,\phi_k\}$ for $\Phi$ satisfies case (ii) of the theorem.

Now suppose that $\psi_1 \neq \phi_{1_L}, \phi'_{1_L}$ and that $D$ is both $\phi_1$- and $\phi'_1$-invariant and $L\cap \fix{\phi_1},
L\cap \fix{\phi'_1}\leq D$. Exactly the same argument as before interchanging $C$ with $D$ and $H$ with $K$ shows that the set of
representatives $\{ \phi_2, \phi_1, \phi_3, \ldots ,\phi_k\}$ for $\Phi$ satisfies case (ii) of the theorem with $H=C$ and
$K=D*\langle z\rangle$.

Next, suppose that $\psi_1 \neq \phi_{1_L}, \phi'_{1_L}$ (say $\psi_1 =\phi_{2_L}$), that $C$ is $\phi_1$-invariant and $L\cap
\fix{\phi_1}\leq C$, and that $D$ is $\phi'_1$-invariant and $L\cap \fix{\phi'_1}\leq D$. In this case let $H=C*D^{z^{-1}}$ and
$y=z\neq 1$. Since $F=L*\langle z\rangle =C*D*\langle z\rangle$, we have $F=H*\langle y\rangle$. The equation
$\phi'_1=\gamma_{z}^{-1}\phi_1 \gamma_z$ and the $\phi'_1$-invariance of $D$ imply that $D^{z^{-1}}$, and so $H$, is
$\phi_1$-invariant. Also, $$\fix{\phi_1}=(L\cap \fix{\phi_1})*(L\cap \fix{\phi'_1})^{z^{-1}}\leq C*D^{z^{-1}}=H.$$ Concerning
$\phi_2$, we know that $C\cap \fix{\phi_2}=H\cap L\cap \fix{\phi_2}=H\cap \fix{\phi_2}$ and
 $$
\begin{array}{rl}
\fix{\phi_2} & =(C\cap \fix{\phi_2})*(D\cap \fix{\phi_2}) \\ & =(H\cap \fix{\phi_2})*\langle d\, \rangle \\ & =(H\cap
\fix{\phi_2})*\langle y^{-1}hy\rangle,
\end{array}
 $$
for some $1\neq d\in D$ and $1\neq h=ydy^{-1}\in H$ (and so, $1\neq H\cap \fix{\phi_2}$). Additionally, $H$ is $\phi_1$-invariant and
intersects $\fix{\phi_2}$ non-trivially so, Lemma~\ref{malnormal} implies that $H$ is also $\phi_2$-invariant. Finally, for each
$i=3,\ldots ,k$, we know that either $C$ is $\phi_i$-invariant and $\fix{\phi_i}\leq C$, or $D$ is $\phi_i$-invariant and
$\fix{\phi_i}\leq D$. In the first case, another application of Lemma~\ref{malnormal} says that $H$ is also $\phi_i$-invariant, and
clearly $\fix{\phi_i}\leq H$. In the second one, we replace $\phi_i$ by the isogredient automorphism $\gamma_z\phi_i \gamma_{z}^{-1}$,
and analogously we deduce that $$\fix{\gamma_z\phi_i \gamma_{z}^{-1}}=(\fix{\phi_i})^{z^{-1}}\leq D^{z^{-1}}\leq H$$ and $H$ is
$(\gamma_z\phi_i \gamma_{z}^{-1})$-invariant. Hence, we conclude that there is a set of representatives for $\Phi$ satisfying case
(iii') of the theorem with $H=C*D^{z^{-1}}$ and $y=z$.

Suppose that $\psi_1 \neq \phi_{1_L}, \phi'_{1_L}$, that $D$ is $\phi_1$-invariant and $L\cap \fix{\phi_1}\leq D$, and that $C$ is
$\phi'_1$-invariant and $L\cap \fix{\phi'_1}\leq C$. Exactly the same argument as before interchanging $\phi_1$ with $\phi'_1$ and
inverting $z$ shows that there is a set of representatives for $\Phi$ satisfying case (iii') of the theorem with $H=C*D^z$ and
$y=z^{-1}$.

So, we have reduced the discussion of the present subcase to the situation where $\psi_1 \in \{\phi_{1_L}, \phi'_{1_L}\}$. We will
only discuss what happens when $\psi_1 =\phi_{1_L}$. For the other possibility, exactly the same arguments interchanging $\phi_1$ with
$\phi'_1$ and inverting $z$ will work. We have a further two possibilities here.

In the first of these possibilities, we consider what happens if $\psi_1 =\phi_{1_L}$, $C$ is $\phi'_1$-invariant and $L\cap
\fix{\phi'_1}\leq C$. Let $H=C*\langle z\rangle$ and $K=D$. As above, $F=H*K$. We already know that $C$ and $D$ are
$\phi_1$-invariant, and hence another application of Lemma~\ref{malnormal} to $\phi_1$ and $\phi'_1=\phi_1 \gamma_g$, says that $g\in
C$ and so $H$ is also $\phi_1$-invariant. Further,
 $$
\begin{array}{rl}
\fix{\phi_1} & =(L\cap \fix{\phi_1})*(L\cap \fix{\phi'_1})^{z^{-1}} \\ & =(C\cap \fix{\phi_1})*(D\cap \fix{\phi_1})* (L\cap
\fix{\phi'_1})^{z^{-1}} \\ & =(H\cap \fix{\phi_1})*(K\cap \fix{\phi_1}),
\end{array}
 $$
where $r(K\cap \fix{\phi_1})=1$, and the last equality is valid because the fact that $(L\cap \fix{\phi'_1})^{z^{-1}}\leq
C^{z^{-1}}\leq H$ justifies one inclusion, while the other inclusion is clear. (Note that the last expression is a free product
because $H\cap K=1$). On the other hand, we know that for each $i=2,\ldots ,k$, either $C$ is $\phi_i$-invariant and
$\fix{\phi_i}=L\cap \fix{\phi_i}\leq C\leq H$, or $D$ is $\phi_i$-invariant and $\fix{\phi_i}=L\cap \fix{\phi_i}\leq D=K$. And in the
first case, another application of Lemma~\ref{malnormal} applied to $\phi_1$ and $\phi_i$ says that $H$ is also $\phi_i$-invariant.
Hence, the set of representatives $\{ \phi_1, \ldots ,\phi_k \}$ for $\Phi$ satisfies case (ii) of the theorem with $H=C*\langle
z\rangle$ and $K=D$.

The second possibility to consider is where $\psi_1 =\phi_{1_L}$, $D$ is $\phi'_1$-invariant and $L\cap \fix{\phi'_1}\leq D$. Let
$H=C*D^{z^{-1}}$ and $y=z\neq 1$. As in the cases above, $F=H*\langle y\rangle$ and $H$ is $\phi_1$-invariant. Also, since $(L\cap
\fix{\phi'_1})^{z^{-1}}\leq D^{z^{-1}}\leq H$,
 $$
\begin{array}{rl}
\fix{\phi_1} & =(L\cap \fix{\phi_1})*(L\cap \fix{\phi'_1})^{z^{-1}} \\ & =(C\cap \fix{\phi_1})*(D\cap \fix{\phi_1})* (L\cap
\fix{\phi'_1})^{z^{-1}} \\ & =(H\cap \fix{\phi_1})*\langle d\, \rangle \\ & =(H\cap \fix{\phi_1})*\langle y^{-1}hy\rangle,
\end{array}
 $$
where $D\cap \fix{\phi_1}=\langle d\, \rangle \neq 1$ and $1\neq h=ydy^{-1}\in H$. On the other hand, for each $i=2,\ldots ,k$, we
know that either $C$ is $\phi_i$-invariant and $\fix{\phi_i}\leq C$, or $D$ is $\phi_i$-invariant and $\fix{\phi_i}\leq D$. In the
first case, another application of Lemma~\ref{malnormal} says that $H$ is also $\phi_i$-invariant, and clearly $\fix{\phi_i}\leq C\leq
H$. In the second one, we replace $\phi_i$ by the isogredient automorphism $\gamma_z\phi_i \gamma_z^{-1}$, and analogously we deduce
that $\fix{\gamma_z \phi_i \gamma_z^{-1}}=(\fix{\phi_i})^{z^{-1}}\leq D^{z^{-1}}\leq H$ and $H$ is $(\gamma_z \phi_i
\gamma_z^{-1})$-invariant. Hence, we conclude that there is a set of representatives for $\Phi$ satisfying case (iii') of the theorem
with $H=C*D^{z^{-1}}$ and $y=z$.

\medskip

\noindent \emph{Subcase 2.3}: writing $\{\phi_{1_L}, \phi'_{1_L}, \phi_{2_L},\ldots ,\phi_{k_L}\}=\{ \psi_1, \ldots ,\psi_{k+1}\}$,
there is a  subgroup $K\leq L$ which is $\psi_j$-invariant for every $j=1,\ldots ,k+1$, and an element $1\neq y\in L$, such that
$L=K*\langle y\rangle$ (so, $r(K)=n-2\geq 1$); moreover, $\fix{\psi_j}\leq K$ if $j\neq 1$, and $\fix{\psi_1}= (K\cap
\fix{\psi_1})*\langle y^{-1}hy\rangle$ for some $1\neq h\in K$.

Let $H=K*\langle z\rangle$. Then, $F=L*\langle z\rangle =K*\langle y\rangle*\langle z\rangle =H*\langle y\rangle$ and we claim that
$H$ is $\phi_i$-invariant for every $i=1,\ldots ,k$. We already know that $K$ is $\phi_i$-invariant. Furthermore, for $i=1$ recall
that $\phi'_1=\gamma_{z}^{-1}\phi_1 \gamma_z =\phi_1 \gamma_g$, where $g=(z\phi_1)^{-1} z\in F$. But $K$ is a free factor of $F$ which
is invariant under $\phi_1$ and $K\cap \fix{\phi'_1}\neq 1$. So, by Lemma~\ref{malnormal}, $g\in K$ and hence $z\phi_1\in H$. Thus,
$H$ is $\phi_1$-invariant. For $i=2,\ldots ,k$, write $\phi_i=\phi_1 \gamma_{x_i}$ for some $x_i\in F$ and the same argument shows
that $x_i\in K$. So, $z\phi_i\in H$ and $H$ is $\phi_i$-invariant.

Suppose that $\psi_1=\phi_{1_L}$. Then, for $i\neq 1$, $\fix{\phi_i}\leq K\leq H$ and we claim that $\fix{\phi_1}=(H\cap
\fix{\phi_1})*\langle y^{-1}hy\rangle$. One of the two inclusions is trivial. For the other, recall that $L\cap \fix{\phi'_1}\leq K$
and so $(L\cap \fix{\phi'_1})^{z^{-1}}\leq H$. Then,
 $$
\begin{array}{rl}
\fix{\phi_1} & =(L\cap \fix{\phi_1})*(L\cap \fix{\phi'_1})^{z^{-1}} \\ & =(K\cap \fix{\phi_1})*\langle y^{-1}hy\rangle *(L\cap
\fix{\phi'_1})^{z^{-1}} \\ & \leq (H\cap \fix{\phi_1})*\langle y^{-1}hy\rangle,
\end{array}
 $$
where $1\neq h\in K\leq H$. Hence, the set of representatives $\{\phi_1, \ldots ,\phi_k\}$ for $\Phi$ satisfies case (iii') of the
theorem.

Suppose that $\psi_1=\phi'_{1_L}$. Exactly the same argument as before interchanging $\phi_{1_L}$ with $\phi'_{1_L}$, and $z$ with
$z^{-1}$, leads to the conclusion that the set of representatives $\{\phi'_1, \phi_2, \ldots ,\phi_k\}$ for $\Phi$ satisfies case
(iii') of the theorem.

Finally, suppose that $\psi_1\neq \phi_{1_L}, \phi'_{1_L}$, say $\psi_1=\phi_{2_L}$. Then, the following subgroups $L\cap
\fix{\phi_1}$, $L\cap \fix{\phi'_1}$, $\fix{\phi_{3}}, \ldots ,\fix{\phi_{k}}$ are all subgroups of $K\leq H$. In particular,
$\fix{\phi_1}=(L\cap \fix{\phi_1})*(L\cap \fix{\phi'_1})^{z^{-1}}\leq H$. And $\fix{\phi_2}=L\cap \fix{\phi_2}=(K\cap
\fix{\phi_2})*\langle y^{-1}hy\rangle =(H\cap \fix{\phi_2})*\langle y^{-1}hy\rangle$, where $1\neq h\in K\leq H$. Hence, the set of
representatives $\{\phi_2, \phi_1, \phi_3, \ldots ,\phi_k\}$ for $\Phi$ satisfies case (iii') of the theorem.

This completes case 2, and the whole proof of the theorem. \qed

\section{The main Theorem}\label{applications}

As immediate consequences of Theorem~\ref{main}, we can deduce Theorem~\ref{mainconnex} concerning single automorphisms, and
Theorem~\ref{cormain} providing a more explicit description of what a 1-auto-fixed subgroup of $F$ looks like.

{\bf \ref{mainconnex}~Theorem }{\it Let $\phi$ be an automorphism of a finitely generated free group $F$. Then, either $\fix{\phi}$ is
cyclic or there exists a non-trivial free factorisation $F=H*K$ such that $H$ is $\phi$-invariant and one of the following holds:
\begin{apartats}
\item $\fix{\phi}\leq H$, \item $K$ is also $\phi$-invariant and
$\fix{\phi}=(H\cap \fix{\phi})*(K\cap \fix{\phi})$, where $r(K\cap \fix{\phi})=1$, \item there exist non-trivial elements $y\in F$,
$h,h'\in H$, such that $K=\langle y\rangle$, $y\phi =h'y$, $h$ is not a proper power, $\fix{\phi}=(H\cap \fix{\phi})*\langle y^{-1}hy
\rangle$ and $h\phi =h'hh'^{-1}$.
\end{apartats}}

\demo If $\fix{\phi}$ is cyclic we are done. Otherwise, apply Theorem~\ref{main} to the outer automorphism $\Phi$ containing $\phi$.
We obtain a set of representatives $\{ \phi_1, \ldots ,\phi_k\}$ for $\Phi$ satisfying~\ref{main}(i), or~\ref{main}(ii)
or~\ref{main}(iii), and such that $\phi =\phi_i$ for some $i$.

If~\ref{main}(i) holds then $\fix{\phi}$ is contained in a proper free factor of $F$. Applying Lemma~\ref{iff}, we end up in case (i).

Suppose that~\ref{main}(ii) holds. If $\phi =\phi_1$ we are in case (ii); otherwise, in case (i).

Finally, suppose~\ref{main}(iii). If $\phi =\phi_1$ we are in case (iii); otherwise, in case (i). \qed

{\bf \ref{cormain}~Theorem } {\it Let $F$ be a non-trivial finitely generated free group and let $\phi\in Aut(F)$ with $\fix{\phi}
\neq 1$. Then, there exist integers $r,s\geq 0$, $\phi$-invariant non-trivial subgroups $K_1, \ldots ,K_r \leq F$, primitive elements
$y_1, \ldots ,y_s \in F$, a subgroup $L\leq F$, and elements $1\neq h'_j \in H_j =K_1*\cdots *K_r*\langle y_1,\ldots ,y_j \rangle$,
$j=0,\ldots ,s-1$, such that
 $$
F=K_1*\cdots *K_r*\langle y_1,\ldots ,y_s \rangle *L
 $$
and $y_j \phi =h'_{j-1}y_j$ for $j=1,\ldots ,s$; moreover,
 $$
\fix{\phi}=\langle w_1, \ldots ,w_r, y_1^{-1}h_0 y_1, \ldots ,y_s^{-1}h_{s-1}y_s \rangle
 $$
for some non-proper powers $1\neq w_i\in K_i$ and some $1\neq h_j\in H_j$ such that $h_j\phi =h_j' h_j h_j'^{-1}$, $i=1,\ldots ,r$,
$j=0,\ldots ,s-1$.}

\demo The proof is by induction on the rank of $F$. If $r(F)=1$, then $\phi$ must be the identity map and the result is clear. So,
suppose $r(F)\geq 2$ and the result known for free groups of smaller rank.

If $\fix{\phi}$ has rank 1, take $r=1$, $K_1=F$, $s=0$ and $L=1$, and we are done.

Otherwise, apply Theorem~\ref{mainconnex} to $\phi$. We obtain a non-trivial free factorisation $F=H*K$ such that $H$ is
$\phi$-invariant and one of~\ref{mainconnex}~(i), or~\ref{mainconnex}~(ii) or~\ref{mainconnex}~(iii) is satisfied. Let us apply the
inductive hypothesis to $\phi_{_H}\in Aut(H)$ (using the notation above) and distinguish the three cases.

If $\fix{\phi}\leq H$ then changing $L$ to $L*K$, we are done.

If $K$ is $\phi$-invariant and $\fix{\phi}=(H\cap \fix{\phi})*(K\cap \fix{\phi})$ with the second factor $K\cap \fix{\phi}=\langle
w\rangle \neq 1$, then increasing $r$ to $r+1$ and adding $K_{r+1}=K\neq 1$ and $w_{r+1}=w$ to the previous structure, we are done.

Finally, if there exist non-trivial elements $y\in F$ and $h,h'\in H$, such that $K=\langle y\rangle$, $y\phi =h'y$, $h$ is not a
proper power, $\fix{\phi}=(H\cap \fix{\phi})*\langle y^{-1}hy \rangle$ and $h\phi =h'hh'^{-1}$, then increasing $s$ to $s+1$ and
adding $y_{s+1}=y$, $h_s'=h'$ and $h_s=h$ to the previous structure, we are done. \qed

Observe that, from the description in Theorem~\ref{cormain} (or an application of Theorem~\ref{mainconnex} and a simple inductive
argument), one can immediately deduce Bestvina-Handel Theorem: $r(\fix{\phi})\leq r(F)$ for every $\phi \in Aut(F)$. In fact, every
$K_i$ contributes with one unity to the left hand side and with $r(K_i)\geq 1$ to the right hand side; every $y_j$ contributes with
one unity to both sides; and $L$ contributes nothing to the left.

\section*{Appendix}

We will justify in this appendix the construction of the automorphism given in \ref{construct}~(iii). First we prove a technical
lemma.

{\bf Lemma} {\it Let $F$ be a free group, $K\leq F$ a finitely generated subgroup, and let $h\neq 1$ and $u_i$, $i\in I$, be elements
of $F$. If $h^{u_i} \in K$ for every $i\in I$ then, the set $\{ u_i K \, \vert \, i\in I\}$ of right cosets of $K$ is finite.}

\demo Pick a basis $X$ for $F$, and consider the coset graph $Z$ of $K\leq F$ with respect to $X$, i.e. the covering with fundamental
group $K$, of the bouquet labelled with the elements of $X$. The vertices of this graph are the right cosets $gK$ of $K$ and, for
every $x\in X$, there is an edge from $gK$ to $xgK$, $g\in F$ (see~\cite{N} for more details). Observe that the finite generation of
$K$ implies that the core $c(Z)$ is finite.

For every $i\in I$, we have $u_i^{-1}hu_i\in K$ and so $hu_iK =u_i K$. Since $h\neq 1$, the $X$-reduced form of $h$ determines a
non-trivial closed path in $Z$ based at $u_iK$. This implies that the distance (in the graph $Z$) from the vertex $u_iK$ to the
subgraph $c(Z)$ is not more than $\vert h\vert_{_X} /2$. The finiteness of $c(Z)$ together with this uniform upper bound imply that
$\{u_iK \,\vert\, i\in I\}$ is a finite set of vertices. \qed

{\bf Proposition} {\it Let $\{ a_1, \ldots ,a_{n-1}\}$ be a basis for $H$, let $\varphi \in Aut(H)$ and suppose that $h\varphi
=h'hh'^{-1}$ for some $1\neq h,h'\in H$ with $h$ not a proper power. By adding a new free generator $y$, we obtain a bigger free group
$F$, and $\varphi$ can be extended to an automorphism $\phi \in Aut(F)$ by setting $a_i\phi_r=a_i\varphi$ for $i=1,\ldots,n-1$, and
$y\phi_r =h'h^r y$. Then, ${\phi_r}_{_{H}}=\varphi$ and, for all but finitely many choices of the integer $r$, the fixed subgroup of
$\phi_r$ is precisely $\fix{\phi_r}=\fix{\varphi}*\langle y^{-1}hy \rangle$. }

\demo  The inclusion $\fix{\phi_r}\geq \fix{\varphi}*\langle y^{-1}hy \rangle$ is clear for every $r$. So, it only remains to show the
other inclusion, for every $r$ except finitely many.

The result holds if $H$ is cyclic as in this case, $h$ is a generator of $H$ and a simple cancellation argument shows the equality is
satisfied for every integer $r$ except that with $h'h^r=1$. Alternatively, it is clear that the only subgroup of $F=\langle h,y
\rangle$ properly containing $\langle h, y^{-1} h y \rangle $ is $F$ itself, which is not equal to  $\fix{\phi_r}$ except when
$h'h^r=1$.

So, assume that $r(H)\geq 2$. Consider the automorphisms of $H$ given by the equation $\varphi_r=\varphi \gamma_{h' h^r}$, $r\in
\mathbb{Z}$. Observe that $h\in \fix{\varphi_r}$ for all $r$, and that $\varphi_r =\varphi_{r'}$ if, and only if, $r=r'$ (since
$r(H)\geq 2$).

Fix a value of $r$ and consider the set $S_r$ of those integers $s$ for which $\varphi_s$ is isogredient to $\varphi_r$. For every
$s\in S_r$, let $u_s \in H$ be such that $\varphi_r =\gamma_{u_s}^{-1}\varphi_s \gamma_{u_s}$. Thus for $s \in S_r$ and every $x \in
H$, $$
\begin{array}{rcl}
(x \varphi)^{h'h^r} & = & x \varphi_r \\ &=& x \gamma_{u_s}^{-1} \varphi_s \gamma_{u_s} \\ &=& x \gamma_{u_s}^{-1} \varphi
\gamma_{h'h^s} \gamma_{u_s}\\ &=& (x \varphi)^{(u_s \varphi)^{-1} h'h^s u_s}. \\
\end{array}
$$ Thus, since $r(H) \geq 2$, we deduce that $u_s\varphi= h'h^s u_s h^{-r} h'^{-1}$. Now, a simple computation shows that, for every
such $s$, $h^{u_s}\in \fix{\varphi_r}$. Since $\fix{\varphi_r}$ is finitely generated, the previous lemma says that the set $\{ u_s
\fix{\varphi_r} \, \vert \, s\in S_r\}$ of right cosets of $\fix{\varphi_r}$ in $H$ is finite. And it is straightforward to verify
that for $s,s'\in S_r$, $u_s\fix{\varphi_r} =u_{s'}\fix{\varphi_r}$ if, and only if, $s=s'$. Hence, $S_r$ is finite, that is, each
$\varphi_r$ is isogredient to only finitely many of the others.

Now, applying Theorem~\ref{bh no connex} to the outer automorphism $\Psi\in Out(H)$ determined by $\varphi_r$, there exist a finite
set $S$ of integers such that $\fix{\varphi_r}$ is cyclic for all $r\not\in S$. Since $h\in \fix{\varphi_r}$ and it is not a proper
power, we have $\fix{\varphi_r}=\langle h\rangle $ for all $r\not\in S$. Moreover, for every such value of $r$, we have $\vert S_r
\vert =1$, since $s\in S_r$ implies $h, h^{u_s} \in \fix{\varphi_r}$ and hence $u_s \in \langle h\rangle $ forcing $\varphi_r
=\varphi_s$ and $r=s$. Thus, we deduce that for $1 \neq r \not\in S$, hence for all but finitely many $r$, $\fix{\varphi_r}=\langle
h\rangle $ and $\varphi_r \neq \gamma_{u}^{-1} \varphi \gamma_u$ for every $u \in H$.

Now, take such a value of $r$ and consider the automorphism $\phi_r\in Aut(F)$ above. A simple cancellation argument shows that any
$\phi_r$-fixed word $w\in F$ can be written as a product of $\phi_r$-fixed words of the form $u$, $(uy)^{\epsilon}$ or $y^{-1}uy$, for
some $u \in H$, $\epsilon=\pm 1$. The second type never occurs since $(uy)\phi_r =uy$ implies $u\varphi =uh^{-r}h'^{-1}$ and
$\gamma_{u}^{-1} \varphi \gamma_u=\varphi_r$, which is not the case. Similarly, one can prove that $(y^{-1}uy)\phi_r =y^{-1}uy$
implies $u\in \fix{\varphi_r}=\langle h \rangle$. Thus, $\fix{\phi_r}=\fix{\varphi}*\langle y^{-1}hy \rangle$, and this is valid for
all but finitely many $r$. \qed

\section*{Acknowledgments}

Both authors thank Warren Dicks for interesting comments and suggestions, that turned out to be the original motivation for this
research. We also thank M. Feighn for very useful comments related to section~\ref{improved rtt}.

The first author gratefully acknowledges support by EPSRC. The second one expresses his gratitude to the City College of New York for
the hospitality received during the academic course 2000-2001, when this paper has been written. He also acknowledges partial support
by the DGI (Spain) through Grant BFM2000-0354, and by the DGR (Generalitat de Catalunya) through grant 2001BEA1400176.

\end{document}